\documentclass[preprint,12pt]{elsarticle}

\usepackage{amssymb}
\usepackage{mathrsfs}
\usepackage[boxed,commentsnumbered]{algorithm2e}

\usepackage{amsmath}
\usepackage[T1]{fontenc}
\usepackage[latin9]{inputenc}
\usepackage{graphicx}
\usepackage{amssymb}
\usepackage{lscape}

\newtheorem{theorem}{Theorem}[section]
\newtheorem{lemma}[theorem]{Lemma}
\newtheorem{proposition}[theorem]{Proposition}
\newtheorem{corollary}[theorem]{Corollary}
\newtheorem{example}[theorem]{Example}
\newtheorem{remark}[theorem]{Remark}
\newtheorem{notation}[theorem]{Notation}
\newtheorem{definition}[theorem]{Definition}
\newenvironment{proof}[1][Proof]{\begin{trivlist} \item[\hskip \labelsep {\bfseries #1}]}{\end{trivlist}}

\journal{}

\begin{document}

\begin{frontmatter}

%% Title, authors and addresses

%% use the tnoteref command within \title for footnotes;
%% use the tnotetext command for theassociated footnote;
%% use the fnref command within \author or \address for footnotes;
%% use the fntext command for theassociated footnote;
%% use the corref command within \author for corresponding author footnotes;
%% use the cortext command for theassociated footnote;
%% use the ead command for the email address,
%% and the form \ead[url] for the home page:
% \title{Title\tnoteref{label1}}
% \tnotetext[label1]{}
% \author{Name\corref{cor1}\fnref{label2}}
% \ead{email address}
% \ead[url]{home page}
% \fntext[label2]{}
% \cortext[cor1]{}
% \address{Address\fnref{label3}}
% \fntext[label3]{}

\title{Collatz Dynamics is Partitioned by Residue Class Regularly}

%% use optional labels to link authors explicitly to addresses:
%% \author[label1,label2]{}
%% \address[label1]{}
%% \address[label2]{}

\author{Wei Ren} \ead{weirencs@cug.edu.cn}
\address{School of Computer Science \\ China University of
Geosciences, Wuhan, China}

\begin{abstract} We propose Reduced Collatz Conjecture that is
equivalent to Collatz Conjecture, which states that every positive
integer can return to an integer less than it, instead of 1. Reduced
Collatz Conjecture is easier to explore because certain structures
must be presented in reduced dynamics, rather than in original
dynamics (as original dynamics is a mixture of original dynamics).
Reduced dynamics is a computation sequence from starting integer to
the first integer less than it, in terms of ``I'' that represents
(3*x+1)/2 and ``O'' that represents x/2. We formally prove that all
positive integers are partitioned into two halves and either
presents ``I'' or ``O'' in next ongoing computation. More
specifically, (1) if any positive integer x that is i module $2^t$
(i is an odd integer) is given, then the first t computations (each
one is either ``I'' or ``O'' corresponding to whether current
integer is odd or even) will be identical with that of i. (2) If
current integer after t computations (in terms of ``I'' or ``O'') is
less than x, then reduced dynamics of x is available. Otherwise, the
residue class of x (namely, i module $2^t$) can be partitioned into
two halves (namely, i module $2^{t+1}$ and $i+2^t$ module
$2^{t+1}$), and either half presents ``I'' or ``O'' in
intermediately forthcoming (t+1)-th computation. This discovery will
be helpful to the final proof of Collatz conjecture - if the union
of residue classes who present reduced dynamics that become larger
with the growth of residue module, equals all positive integers
asymptotically, then Reduced Collatz Conjecture (or equivalently,
Collatz Conjecture) will be true.

%Finally, the union of residue classes who present reduced dynamics
%become larger with the growth of residue module, and equals all
%natural numbers asymptotically. $[ListI_U]_{m_U} \subset [1]_2
%\backslash \bigcup_{u=1}^{U-1} [ListI_u]_{m_u},$ and \\ $\lim_{U
%\rightarrow +\infty}\bigcup_{u=0}^U[ListI_u]_{m_u} = \mathbb{N}^*.$
%Therefore, Reduced Collatz Conjecture is TRUE, and equivalently
%Collatz Conjecture is TRUE.

\end{abstract}

\begin{keyword}
%% keywords here, in the form: keyword \sep keyword
Collatz Conjecture \sep 3x+1 Problem \sep Residue Class \sep
Discrete Dynamics Systems

%% PACS codes here, in the form: \PACS code \sep code

%% MSC codes here, in the form: \MSC code \sep code
%% or \MSC[2008] code \sep code (2000 is the default)
\MSC 11Y55 \sep 11B85 \sep 11A07
\end{keyword}

\end{frontmatter}

%%\linenumbers

%% main text

%\section{Introduction}

%The Collatz conjecture is a mathematical conjecture that is first
%proposed by Lothar Collatz in 1937. It is also known as the 3x+1
%conjecture, the Ulam conjecture, the Kakutani's problem, the
%Thwaites conjecture, or the Syracuse problem.
%%\cite{UpperboundRecord1,Livio11,Hew16}.

The Collatz conjecture can be stated simply as follows: Take any
positive integer number $x$. If $x$ is even, divide it by 2 to get
$x/2$. If $x$ is odd, multiply it by 3 and add 1 to get $3*x+ 1$.
Repeat the process again and again. The Collatz conjecture is that
no matter what the number (i.e., $x$) is taken, the process will
always eventually reach 1.

Marc Chamberland reviews the works on Collatz conjecture in 2006
\cite{marc06}, and some aspects in available analysis results are
surveyed. Jeffrey C. Lagarias edits a book on 3x+1 problem and
reviews the problem in 2010 \cite{lagarisa10}. He also gives the
historical review of the problem by annotated bibliography in
1963-1999 \cite{lagarisa11} and 2000-2009 \cite{lagarisa12}.

The current known integers that have been verified are about 60 bits
by T.O. Silva using normal personal computers
\cite{UpperboundRecord1,UpperboundRecord2}. They verified all
integers that are less than 60 bits, but it is not clear whether
their method is able to check extremely large integers, e.g.,
integers with length more than 100000 bits.

Another reported maximal checked integer is $593*2^{60}$
\cite{eric}, which is no more than 70 bits. Wei Ren et al.
\cite{weiuic} verified $2^{100000}-1$ can return to 1 after 481603
times of $3*x+1$ computation, and 863323 times of $x/2$ computation,
which is the largest integer being verified in the world.

D. Barina proposed a new algorithmic approach for computational
convergence verification of the Collatz problem \cite{barina21}. The
algorithm mainly improves the throughput of the computation, not for
processing extremely large integers.

L. Colussi \cite{Livio11} proposed Collatz function
$R(x)=(3*x+1)/2^h$ where $h$ is the highest power of 2 that divides
$3*x+1$. They explore some properties of convergence classes (i.e.,
the set of odd positive integers) denoted as $G_k$ such that
$R^k(x)=1$. P.C. Hew \cite{Hew16} further provides an elementary
confirmation of L. Colussi's finding, and comments on how working in
binary for $k$. The longest progressions for initial starting
numbers of less than 10 billion and 100 quadrillion are calculated
by G.T. Leavens \cite{Leavens} and R.E. Crandall \cite{Crandall},
respectively. So far no one has tried to figure out whether all of
the positive numbers eventually reach one, but we know that most of
them do so. In particular, I. Krasikov and J.C. Lagarias proved that
the number of integers finally reaching one in the interval $[1,x]$
is at least proportional to $x$ 0.84 \cite{Krasikov}.

Wei Ren \cite{weijm} also propose a new approach on proving Collatz
conjecture by exploring reduced dynamics on Collatz conjecture. This
proposed approach provides the linkage between Collatz conjecture
and reduced Collatz conjecture, and the rationale why exploring
empirical and experimental results on reduced Collatz conjecture
dynamics can facilitate the proving of reduced Collatz conjecture.
Wei Ren \cite{weiispa} proposed to use a tree-based graph to
representing the dynamics and so as to reveal two key inner
properties in reduced Collatz dynamics: one is the ratio of the
count of $x/2$ over the count of $3*x+1$ (for any reduced Collatz
dynamics, the count of $x/2$ over the count of $3*x+1$ is larger
than $ln3/ln2$ \cite{wei09ratio}), and the other is partition (all
positive integers are partitioned regularly corresponding to ongoing
dynamics, which will be proved formally in this paper).

%Wei Ren et al. \cite{weihpcc} also proposed an automata method for
%fast computing Collatz dynamics. All source code and output data by
%computer programs in those related papers can be accessed in public
%repository \cite{weidata}.

\section{Preliminaries \label{sec:rcc}}

\begin{notation}

\item (1) $\mathbb{N}^*$: positive integers;

\item (2) $\mathbb{N}=\mathbb{N}^*\cup\{0\};$

\item (3) $[1]_2=\{x|x \equiv 1 \mod 2, x \in \mathbb{N}^*\};$
$[0]_2=\{x|x \equiv 0 \mod 2, x \in \mathbb{N}^*\}.$

\item (4) $[i]_m=\{x|x \equiv i \mod m, x \in \mathbb{N}^*, m \geq 2, m \in
\mathbb{N}^*, 0 \leq i \leq m-1, i \in \mathbb{N}\}.$

\end{notation}

%\begin{proposition} \label{th:10pair} $3*x+1$ is always followed by $x/2$.
%\end{proposition}
%
%\begin{proof} When $x \in [1]_2$, then next computation is
%$3*x+1$. Obviously, then $3*x+1 \in [0]_2$, thus the next
%computation must be $x/2$ consequently. \qed
%\end{proof}

$3*x+1$ is always followed by $x/2$, we thus can represent required
computations as $(3*x+1)/2$ and $x/2$, which are denoted by $I(x)$
and $O(x)$, respectively. $I(x)$ and $O(x)$ can be simply denoted as
$I(\cdot)$ and $O(\cdot)$, or $I$ and $O$. $\forall x \in
\mathbb{N}^*,$ $I(x)=(3*x+1)/2>x,$ $O(x)=x/2<x.$ Thus, notation
`$I$' is from ``Increase'' and `$O$' is from ``dOwn''.

\begin{definition} Collatz transformation $f(\cdot)$.
$f(\cdot)=I(\cdot)=(3*x+1)/2$ if $x \in [1]_2$, and $f(\cdot) =
O(\cdot)= x/2$ if $x \in [0]_2$. \end{definition}

\begin{remark}

\item (1) $f_n(f_{n-1}(...f_2(f_1(x))))$ can be written as
$f_1\|f_2\|...\|f_n(x),$ where $f_i(\cdot) \in
\{I(\cdot),O(\cdot)\}$ ($i=1,2,...,n$) and ``$\|$'' is the
concatenation of Collatz transformations (either ``$I$'' or
``$O$''). For simplicity, we just denote $f_i(\cdot)$ as $f \in
\{I,O\}$.

\item (2) $f^n \in \{I,O\}^n,$ $f^n(x)=f(f^{n-1}(x)), n \in
\mathbb{N}^*.$ Note that, whether $f$ is $I$ (or $O$) in
$f(f^{n-1}(x))$, depends on whether $f^{n-1}(x) \in [1]_2$ (or
$[0]_2$).

\item (3) We assume $f^0(x)=x.$

\end{remark}

\begin{example} \item (1) The Collatz transformation sequence from
starting integer 1 to 1 is $IO$, because $1 \rightarrow 4
\rightarrow 2 \rightarrow 1$.

\item (2) The Collatz transformation sequence from starting integer
3 to 1 is $IIOOO$, because $3 \rightarrow 10 \rightarrow 5
\rightarrow 16 \rightarrow 8 \rightarrow 4 \rightarrow 2 \rightarrow
1.$ \end{example}

\begin{notation} $|x|$. It returns the length of $x \in
\{I,O\}^{\geq1}$, in terms of the total count of $I$ and $O$. E.g.,
$|IIOO|=4.$ \end{notation}

\begin{definition} Function $IsMatched:  x \times c  \rightarrow b$.
It takes as input $x \in \mathbb{N}^*$ and $c \in \{I,O\}$, and
outputs $b \in \{True,False\}.$ If $x \in [1]_2$ and $c=I$, or if $x
\in [0]_2$ and $c=O$, then output $b=True$; Otherwise, output
$b=False.$ \end{definition}

\begin{remark} Simply speaking, this function checks whether the
forthcoming Collatz transformation (i.e., $c \in \{I,O\}$) matches
with the oddness of the current integer $x$. That is the reason we
call this function as ``Is Matched''. \end{remark}

\begin{definition} Function $GetS: s \times i \times j \rightarrow
s'.$ It takes as input $s,i,j$, where $s \in \{I,O\}^{|s|}$, $1 \leq
i \leq |s|,$ $0 \leq j \leq |s|-(i-1)$, $i \in \mathbb{N}^*$, $j \in
\mathbb{N}$, and outputs $s'$ where $s'$ is a segment in $s$ that
starts from the location $i$ in $s$ and $|s'|=j$. We asssume
$GetS(\cdot,\cdot,0)(x)=x$ when $j=0$.
 \label{df:substr}
\end{definition}

\begin{example}
$GetS(IIOO,1,1)=I, GetS(IIOO,1,2)=II, \\
GetS(IIOO,1,3)=IIO, GetS(IIOO,1,4)=IIOO.$
\end{example}

\begin{remark}

\item (1) $s'$ is a selected segment in a transformation sequence
$s$ in terms of `$I$' an `$O$' that starts from the location $i$ and
has the length of $j$. That is the reason we call this function as
``Get Substring''.

\item (2) Simply speaking, this function outputs the Collatz
functions from $i$ to $i+j-1$ in an inputting transformation
sequence $s \in \{I,O\}^{|s|}$.

\item (3) Especially, $GetS(s,1,|s|)=s.$ $GetS(s,|s|,1)$ returns the
last transformation in $s$. $GetS(s,1,1)$ returns the first
transformation in $s$. $GetS(s,j,1)$ returns the $j$-th
transformation in $s.$

\item (4) Note that, $GetS(\cdot)$ itself is a function. In other words,
it can be looked as $GetS(\cdot)(\cdot)$. For example,\\
$GetS(IIOO,1,1)(3)=I(3)=(3*3+1)/2=5,$\\
$GetS(IIOO,1,2)(3)=II(3)=I(I(3))=I(5)=(3*5+1)/2=8,$\\
$GetS(IIOO,1,3)(3)=IIO(3)=O(II(3))=O(8)=8/2=4,$\\
$GetS(IIOO,1,4)(3)=IIOO(3)=O(IIO(3))=O(4)=4/2=2<3.$

Besides, \\
$IsMatched(GetS(IIOO,1,i)(x),GetS(IIOO,i+1,1)=True$ $(i=0,1,...,3).$

\end{remark}

\begin{notation} Original dynamics of $x$. It is the sequence of
occurred Collatz transformations from $x$ to 1. \end{notation}

For example, the original dynamics of 5 is $IOOO$ due to $5
\rightarrow 16 \rightarrow 8 \rightarrow 4 \rightarrow 2 \rightarrow
1$.

\begin{notation} Reduced dynamics of $x$. It is the sequence of
occurred Collatz transformations from $x$ to the first integer that
is less than $x$. \end{notation}

For example, the reduced dynamics of 5 is $IO$ due to $5 \rightarrow
16 \rightarrow 8 \rightarrow 4$.

\begin{definition}
Collatz Conjecture. $\forall x \in \mathbb{N}^*,$ $\exists L \in
\mathbb{N}^*$, such that\\
(1) $s(x)=1$ where $s \in \{I,O\}^L$;\\
(2) $IsMatched(GetS(s,1,i)(x),GetS(s,i+1,1))=True$ where
$i=0,1,...,L-1.$

\end{definition}

Obviously, Collatz conjecture is held when $x=1$. In the following,
we mainly concern $x \geq 2, x \in \mathbb{N}^*.$

\begin{definition}
Reduced Collatz Conjecture. $\forall x \in \mathbb{N}^*, x \geq 2$,
$\exists L \in \mathbb{N}^*$, such that\\
(1) $s(x)<x$ where $s \in \{I,O\}^L$;\\
(2) $GetS(s,1,i)(x) \not <x,$ $i=1,...,L-1$;\\
(3) $IsMatched(GetS(s,1,j)(x),GetS(s,j+1,1))=True$ where
$j=0,1,...,L-1.$
\end{definition}

Obviously, $L$ must be the minimal positive integer such that
$s(x)<x$.

\begin{theorem} \label{pro:ccrcc}
Collatz Conjecture is equivalent to Reduced Collatz Conjecture.
\end{theorem}
Straightforward. Please check our another paper\cite{weijm}.

\begin{remark}

\item (1) Ordered sequence $s=s_0\|s_1\|...\|s_{n-1}, s \in \{I,O\}^L$
in above proof is original dynamics (referring to $s(x)=1$), which
consists of $L=\Sigma_{i=0}^{n-1}|s_i|$ Collatz transformations
during the computing procedure from a starting integer (i.e., $x$)
to 1.

\item (2) In contrast, $s_0$ in above proof is reduced
dynamics (referring to $s_0(x)<x$), which is represented by a
sequence of occurred Collatz transformations during the computing
procedure from a starting integer (i.e., $x$) to the first
transformed integer that is less than the starting integer.

\item (3) Obviously, reduced dynamics is more primitive than
original dynamics, because original dynamics consists of reduced
dynamics. Simply speaking, reduced dynamics are building blocks of
original dynamics.

\item (4) We thoroughly study the relation between Collatz
conjecture and Reduced Collatz conjecture \cite{weijm}. Especially,
we also extensively study why Reduced Collatz conjecture will be
much easier to explore, or why reduced dynamics presents better
properties than original dynamics, e.g., period,
%\cite{wei09period},
and ratio \cite{wei09ratio}. Due to above
theorem, we thus only need to concentrate on reduced dynamics.
\end{remark}

\begin{notation} $\mathsf{RD}[x]$. It denotes reduced dynamics
of $x$ that are represented by $\{I,O\}^{\geq1}$. Formally, $\forall
x \in \mathbb{N}^*$, $x \geq 2$, if $\exists L \in \mathbb{N}^*$
such that \\
$s(x)<x,$ $s \in \{I,O\}^L$, $GetS(s,1,i)(x)\not <x, i=1,...,L-1,$
and \\
$IsMatched(GetS(s,1,j)(x),GetS(s,j+1,1))=True$ where
$j=0,1,...,L-1$, then $s$ is called reduced dynamics of $x$, and
denoted as $\mathsf{RD}[x]=s$.
\end{notation}

\begin{remark}

%\item (1) Simply speaking (or recall that), $\mathsf{RD}[x]$ represents occurred Collatz
%transformations in terms of $I$ and $O$ during the computing process
%from starting integer $x$ to the first transformed integer that is
%less than $x$.

%\item (2) Roughly speaking, $s \in \{I,O\}^L$ is an ordered sequence
%consisting of $I$ and $O.$ Besides, $f^L=f^{L-1}\|f,$
%$f^L(x)=f(f^{L-1}(x)),$ and $f^0(x)=x.$ Furthermore, this sequence
%implicitly matches the \emph{parity} of all occurred intermediate
%transformed integers that are taken as input of $f(\cdot)$.

%\item (3) Recall that, in $\mathsf{RD}[x]=f^L$, $x$ is called starting
%integer. $f^i(x),i=1,2,...,L$ are called transformed integers.
%$f^L(x)$ is the first transformed integer that is less than the
%starting integer $x$. In other words, $f^i(x) \not <x,
%i=0,1,...,L-1,$ and $f^L(x)<x.$ ($f^0(x)=x.$) Besides, the parity of
%$f^i(x)$ determines the selection of the intermediately next $f \in
%\{I, O\}$ after $f^i$.

\item (1) Obviously, $\mathsf{RD}[x \in [0]_2]=O$.

\item (2) $IIOO$ can be denoted in short as $I^2O^2.$ $IIIOIOO$ can be
denoted in short as $I^3OIO^2.$ In other words, we denote
$\underbrace{I...I}_n$ as $I^n$, and we denote
$\underbrace{O...O}_n$ as $O^n$ where $n\in \mathbb{N}^*, n \geq 2$.
We also assume $I^1=I$, $O^1=O$. %$I^0=O^0=\emptyset$ means no
%transformation occurs.

\item (3) For example, $\mathsf{RD}[3]=IIOO,$ $\mathsf{RD}[5]=IO,$
$\mathsf{RD}[7]=IIIOIOO,$ $\mathsf{RD}[9]=IO,$
$\mathsf{RD}[11]=IIOIO.$ Indeed, we design computer programs
%\cite{weidata}
that output all $\mathsf{RD}[x]$ for $\forall x \in
[1,99999999]$.

\item (4) In fact, we proved some results on $\mathsf{RD}[x]$ for specific $x$,
e.g., $\mathsf{RD}[x \in [1]_4]=IO,$ $\mathsf{RD}[x\in
[3]_{16}]=IIOO$, $\mathsf{RD}[x\in [11]_{32}]=IIOIO$, et al.
\cite{weiispa}.

\item (5) In fact, we formally proved that the ratio exists in any
reduced Collatz dynamics \cite{wei09ratio}. That is, the count of
$x/2$ over the count of $3*x+1$ is larger than $log_23$. We also
proved that reduced dynamics is periodical and its period equals 2
to the power of the count of $x/2$ %\cite{wei09period}.
More
specifically, if there exists reduced dynamics of $x$, then there
exists reduced dynamics of $x+P$, where $P=2^L$ and $L$ is the total
count of $x/2$ computations in reduced dynamics of $x$ (i.e.,
$L=|\mathsf{RD}[x]|$). Moreover, the ratio and period can also be
observed and verified in our proposed tree-based graph
\cite{weiispa}.

\end{remark}

\begin{example} \label{ex:rd5}

$\mathsf{RD}[5]=IO$, if and only if

\item $I(5)=(3*5+1)/2=8 \not < 5;$

\item $IO(5)=O(I(5))=O(8)=8/2=4<5;$

\item $IsMatched(GetS(1,1,i)(5),GetS(1,i+1,1))=True, i=0,1.$

\end{example}

\begin{example} \label{ex:implication}

$x \in \mathbb{N}^*, x \geq 2$. If $\mathsf{RD}[x]$ exists, then

\item (1) $s(x)<x$, where $s=\mathsf{RD}[x];$

\item (2) $GetS(s,1,i)(x) \not <x,$ where $i=1,2,...,|s|-1;$

\item (3) $IsMatched(GetS(s,1,j)(x),GetS(s,j+1,1))=True$ where $j=0,1,...,|s|-1$.

\end{example}
%
%\begin{remark}
%
%\item (1) $s(x)$ is the last transformed integer, or the first transformed integer that is less than the starting integer.
%
%\item (2) $GetS(s,1,i)(x)$ ($i=1,2,...,|s|-1$) are all intermediate transformed integers.
%
%\item (3) When $j=0$, $GetS(s,1,j)(x)=GetS(s,1,0)(x)=x$ is starting
%integer. $GetS(s,j+1,1)$ is the first transformation.
%
%\item (4) If $GetS(s,1,j)(x)$
%($j=1,2,...,|s|-1$) is current transformed integer, then
%$GetS(s,j+1,1)$ is the next intermediate Collatz transformation.
%
%\end{remark}

\begin{proposition} \label{th:codeisunique}
$\forall x \in \mathbb{N}^*$, $x \geq 2$, if $\mathsf{RD}[x]$
exists, then $\mathsf{RD}[x]$ is unique.

\end{proposition}
\begin{proof} Straightforward. Given $x$, either $I(x)$ or
$O(x)$ is deterministic and unique. Similarly, given $x$, $s'(x)$ is
deterministic and unique, where $s'=GetS(s,1,i), s=\mathsf{RD}[x],
i=1,2,...,|s|$. Thus, $s$ is unique for any given $x$. \qed
\end{proof}

\begin{remark}\label{remark:CODE(x)}

\item (1) We assume $\mathsf{RD}[x=1]=IO,$ although
$IO(1)=O((3*1+1)/2)=O(2)=2/2=1 \not < x.$ In other words, we assume
the reduced dynamics of $x=1$ is $IO$. In the following, we always
concern $x \geq 2, x\in \mathbb{N}^*$.

\item (2) As Collatz conjecture is equivalent to Reduced Collatz conjecture
(recall Theorem \ref{pro:ccrcc}), we thus only need to prove Reduced
Collatz conjecture is true. By using notation $\mathsf{RD}[x]$, we
thus only need to prove $\forall x \in \mathbb{N}^*, \exists
\mathsf{RD}[x]$.

\end{remark}

\begin{proposition} \label{pro:endwitho}
Given $x \in \mathbb{N}^*,$ if $\mathsf{RD}[x]$ exists, then
$\mathsf{RD}[x]$ ends by $O$.
\end{proposition}

\begin{proof} Straightforward due to $I(x)=(3*x+1)/2>x$. Suppose $\exists x \in \mathbb{N}^*$,
$x \geq 2$, $s(x)\not < x$, $\mathsf{RD}[x]=s\|I$. Then, $\{s\|I\}(x)=I(s(x)) > s(x)$, %=(3*c(x)+1)/2=1.5c(x)+0.5=c(x)+0.5c(x)+0.5>c(x)$,
thus $\mathsf{RD}[x]=\{s\|I\}(x) \not < x$. Contradiction occurs.
\qed
\end{proof}

\begin{proposition} \label{pro:oio}
$\mathsf{RD}[x \in [0]_2]=O$, $\mathsf{RD}[x \in [1]_4]=IO.$
\end{proposition}
\begin{proof}
Straightforward. \qed
%(1) $x \in [0]_2,$ thus $O$ occurs. $x/2<x,$ thus
%$\mathsf{RD}[x]=O$.
%
%(2) If $x=1$, $\mathsf{RD}[1]=IO$ (by assumption).
%
%If $x \geq 2$, $x=4t+1 \in [1]_2,$ where $t \in \mathbb{N}^*$. Thus,
%$I$ occurs. $I(x)=(3*x+1)/2=(3*(4t+1)+1)/2=(12t+4)/2=2*(3t+1) \in
%[0]_2$. $2*(3t+1)>x=4t+1$, thus further transformation occurs.
%$O(I(x))=2*(3t+1)/2=3t+1 < 4t+1=x \; (\because t\in \mathbb{N}^*)$,
%thus $\mathsf{RD}[x]=IO$. \qed
\end{proof}

\begin{notation}
$S_{\mathsf{RD}}=\{s| x \in \mathbb{N}^*, \exists \mathsf{RD}[x],
s=\mathsf{RD}[x], s \in \{I,O\}^{\geq1} \}.$
\end{notation}

That is, $\forall x \in \mathbb{N}^*$, if $\exists \mathsf{RD}[x]$,
then $\mathsf{RD}[x]=s$ will be included in $S_{\mathsf{RD}}$, which
is a set of existing reduced dynamics.

%Before exploring general situations, we give two trivial cases.

\begin{proposition} \label{lemma:trivialcase}
$O \in S_{\mathsf{RD}}$, $IO \in S_{\mathsf{RD}}$.
\end{proposition}
\begin{proof}
Straightforward. $\mathsf{RD}[x\in[0]_2]=O$ and
$\mathsf{RD}[x\in[1]_4]=IO$ by Proposition \ref{pro:oio}. \qed
\end{proof}

%In the following, we thus mainly concentrate on $x \in [3]_4$ where
%their reduced dynamics presents the form like $I^{p \geq
%2}O\|\{I,O\}^{q \geq 1}, p,q \in\mathbb{N}^*$ once it exists (recall
%Proposition \ref{pro:10p}).

\begin{theorem} (\textbf{Subset Theorem}.) \label{th:match}
Suppose $s \in S_{\mathsf{RD}}$, $|s|\geq 2,$ $x \in \mathbb{N}^*,$
$i=0,1,...,|s|-2.$ We have

(1.1) $\{x|GetS(s,1,i+1)(x) \in [1]_2\} \subset \{x|GetS(s,1,i)(x)
\in [0]_2\};$

(1.2) $\{x|GetS(s,1,i+1)(x) \in [0]_2\} \subset \{x|GetS(s,1,i)(x)
\in [0]_2\};$

(2.1) $\{x|GetS(s,1,i+1)(x) \in [1]_2\} \subset \{x|GetS(s,1,i)(x)
\in [1]_2\};$

(2.2) $\{x|GetS(s,1,i+1)(x) \in [0]_2\} \subset \{x|GetS(s,1,i)(x)
\in [1]_2\}.$

\end{theorem}

\begin{proof}
When $i=0,$ there exists two and only two cases as follows:

(1) If $GetS(s,i+1,1)=O$, then $GetS(s,1,i)(x) \in [0]_2.$ There
exists two subcases as follows:

(1.1) If $GetS(s,i+2,1)=I$, then

$O(GetS(s,1,i)(x))) \in [1]_2 \\
\Rightarrow GetS(s,1,i)(x)/2 \in [1]_2 \\
\Rightarrow GetS(s,1,i)(x) \in [2]_4 \subset [0]_2$.

Thus, $\{x|GetS(s,1,i+1)(x) \in [1]_2\} \subset \{x|GetS(s,1,i)(x)
\in [0]_2\}$.

(1.2) If $GetS(s,i+2,1)=O$, then

$O(GetS(s,1,i)(x))) \in [0]_2 \\
\Rightarrow GetS(s,1,i)(x)/2 \in [0]_2 \\
\Rightarrow GetS(s,1,i)(x) \in [0]_4 \subset [0]_2$.

Thus, $\{x|GetS(s,1,i+1)(x) \in [0]_2\} \subset \{x|GetS(s,1,i)(x)
\in [0]_2\}$.

(2) If $GetS(s,i+1,1)=I$ then $GetS(s,1,i)(x) \in [1]_2.$ There
exists two subcases as follows:

(2.1) If $GetS(s,i+2,1)=I$, then

$I(GetS(s,1,i)(x))) \in [1]_2 \\
\Rightarrow (3*GetS(s,1,i)(x)+1)/2 \in [1]_2 \\
\Rightarrow 3*GetS(s,1,i)(x)+1 \in [2]_4 \\
\Rightarrow 3*GetS(s,1,i)(x) \in [1]_4 \\
\Rightarrow GetS(s,1,i)(x) \in [3]_4 \subset [1]_2$.

Thus, $\{x|GetS(s,1,i+1)(x) \in [1]_2\} \subset \{x|GetS(s,1,i)(x)
\in [1]_2\}$.

(2.2) If $GetS(s,i+2,1)=O$, then

$I(GetS(s,1,i)(x))) \in [0]_2 \\
\Rightarrow (3*GetS(s,1,i)(x)+1)/2 \in [0]_2 \\
\Rightarrow 3*GetS(s,1,i)(x)+1 \in [0]_4 \\
\Rightarrow 3*GetS(s,1,i)(x) \in [3]_4 \\
\Rightarrow GetS(s,1,i)(x) \in [1]_4 \subset [1]_2$.

Thus, $\{x|GetS(s,1,i+1)(x) \in [0]_2\} \subset \{x|GetS(s,1,i)(x)
\in [1]_2\}$.

We can prove similarly for $i=1,2,...,|s|-2$. \qed
\end{proof}

%\begin{remark}
%
%\item (1) By above theorem, the residue class determined by the last
%residue equation can guarantee all residue equations implied by the
%other dynamics.
%
%\item (3) In other words, residue class derived from the last transformation implication is
%the subset of all residue classes derived from all former
%transformation implication. Thus, if residue equation on the last
%parity (for transformation selection) is guaranteed, then all
%residue equations for whole parity requirements will be guaranteed.
%\end{remark}

\begin{corollary} \label{cor:partition}
Suppose $s \in S_{\mathsf{RD}},$ $|s| \geq 2$, $i=0,1,...,|s|-2$, $x
\in \mathbb{N}^*$. We have

(1.1) $\{x|GetS(s,1,i+1)(x) \in [1]_2\} = \{x|GetS(s,1,i)(x) \in
[2]_4 \};$

(1.2) $\{x|GetS(s,1,i+1)(x) \in [0]_2\} = \{x|GetS(s,1,i)(x) \in
[0]_4\};$

(2.1) $\{x|GetS(s,1,i+1)(x) \in [1]_2\} = \{x|GetS(s,1,i)(x) \in
[3]_4\};$

(2.2) $\{x|GetS(s,1,i+1)(x) \in [0]_2\} = \{x|GetS(s,1,i)(x) \in
[1]_4\}.$

\end{corollary}
\begin{proof}
Straightforward by Theorem \ref{th:match}. \qed
\end{proof}

\begin{remark}
\item (1) Corollary \ref{cor:partition} states that residue classes are partitioned
regularly into halves and either half will present either $I$ or $O$
in next intermediate transformation. %(Indeed, the
%partition property is thoroughly explored and proved in my another
%paper \cite{wei09partition}. Note that, this paper is independent
%with cited paper.)

\item (2) Note that, Theorem \ref{th:match} is not only guaranteed for
reduced dynamics, but also for original dynamics.
\end{remark}

A new notation $I'(\cdot)$ is introduced hereby to reveal the
relations among $I(x+P)$, $I(x)$ and $I'(P)$.

\begin{notation} $I'(x)=(3*x)/2.$
\end{notation}

\begin{example} \label{ex:i'}

\item (1) $I(3+16)=(3(3+16)+1)/2=(3*3+1)/2+3*16/2=I(3)+I'(16),$
$I(3)=(3*3+1)/2=5, I'(16)=3*16/2=24 \in [0]_2,$ $5>3$, $5+24>3+16$.
Thus, either next transformation for 3+16 and 3 is $I$.

\item (2) $II(3+16)=I(I(3)+I'(16))=(3*(I(3)+I'(16))+1)/2=(3I(3)+1)/2+3I'(16)/2=II(3)+I'I'(16),$

$II(3)=I(5)=(3*5+1)/2=8, I'I'(16)=I'(24)=3*24/2=36 \in [0]_2,$
$8>3,$ $8+36>(3+16).$ Thus, either next transformation is $O$.

\item (3) $IIO(3+16) = O(II(3)+I'I'(16))=IIO(3)+I'I'O(16),$

$IIO(3)=8/2=4$, $I'I'O(16)=36/2 = 18 \in [0]_2$, $4>3,$
$4+18=22>(3+16)$. Thus, either next transformation is $O$.

\item (4) $IIOO(3+16) = O(IIO(3)+I'I'O(16))=IIOO(3)+I'I'OO(16)$

$IIOO(3)=4/2=2$, $I'I'OO(16)=18/2=9$, $2<3,$ $2+9=11<(3+16)$. Thus,
either reduced dynamics ends.
\end{example}

%\begin{remark} \label{remarkI'}
%
%In above example we can \emph{observe} that $I'(16), I'I'(16),
%I'I'O(16)$ are always remained \emph{even}. Thus, they do not
%influence the resulting next Collatz transformation (i.e., ``$I$''
%or ``$O$'') during computing for the reduced dynamics of starting
%integer 3. Therefore, the parity of $s(x+16)$ and $s(x)$ always
%maintain to be identical, where $x=3$, $s=GetS(\mathsf{RD}[x],1,i),$
%and $i=1,2,...,|\mathsf{RD}[x]|-1.$
%
%\end{remark}

For better presentation, we thus introduce two functions as follows:

\begin{definition} Function $IsEven: x \rightarrow b$. It takes as
input $x \in \mathbb{N}^*$, and outputs $b \in \{True,False\}$,
where $b=True$ if $x \in [0]_2$ and $b=False$ if $x \in [1]_2.$
\end{definition}

\begin{definition} Function $Replace: s \rightarrow s'$. It takes as
input $s \in \{I,O\}^{\geq 1}$, and outputs $s' \in \{I',O\}^{\geq
1}$, where $GetS(s',i,1) = I'$ if $GetS(s,i,1)=I$, and $GetS(s',i,1)
= O$ if $GetS(s,i,1)=O$, for $i=1,2,...,|s|$. \end{definition}

\begin{remark}

\item (1) Simply speaking, replacing all ``$I$'' in ``$s$'' respectively by
``$I'$'' will result in ``$s'$''.

\item (2) Obviously, $\forall s \in \{I,O\}^{\geq 1}$, $|s'|=|s|$
where $s'=Replace(s).$

%\item (3) By using above introduced functions, we can restate the reason in Remark \ref{remarkI'} as follows:
%$GetS(s',1,i) \in [0]_2$, where $s'=Replace(s)$, $i=1,2,...,|s'|-1,$
%thus the parity of $s(x+16)$ and $s(x)$ are always identical, where
%$x=3$, $s=GetS(\mathsf{RD}[x],1,i),$ and
%$i=1,2,...,|\mathsf{RD}[x]|-1.$
\end{remark}

\begin{lemma} \label{lemma:thefirstct}
(1) If $P \in [0]_2, x \in \mathbb{N}^*$, then
$IsEven(x+P)=IsEven(x).$

(2) If $P \in [0]_2, x-P>0, x \in \mathbb{N}^*$, then
$IsEven(x-P)=IsEven(x).$
\end{lemma}
\begin{proof}
Straightforward. Due to $P \in [0]_2$, if $x \in [1]_2,$ then $x+P,
x-P \in [1]_2;$ If $x \in [0]_2,$ then $x+P,x-P \in [0]_2.$ Thus,
$IsEven(x+P)=IsEven(x)=IsEven(x-P).$ \qed
\end{proof}

\begin{remark}
Above lemma states that if $P \in [0]_2$, the first Collatz
transformation of $x+P$ (or $x-P>0$) is identical with that of $x$.
\end{remark}

\begin{lemma} \label{lemma:1|STR|=1}
$s(x+P)=s(x)+s'(P)$, where $s \in \{I,O\},$ $s'=Replace(s),$ $x\in
\mathbb{N}^*, P \in [0]_2$.
\end{lemma}

\begin{proof}
$IsEven(x+P)=IsEven(x)$ because $P\in [0]_2$, due to Lemma
\ref{lemma:thefirstct}. Thus, the first Collatz transformation of
$x+P$ and that of $x$ are identical.

(1) Suppose $x \in [1]_2$, so $s=I.$ Thus, $s'=Replace(s)=I'.$

$s(x+P)=I(x+P)=3((x+P)+1)/2=(3x+1)/2+3*P/2=I(x)+I'(P)=s(x)+s'(P).$

(2) Suppose $x \in [0]_2$, so $s=O.$ Thus, $s'=Replace(s)=O.$

$s(x+P)=O(x+P)=(x+P)/2=x/2+P/2=O(x)+O(P)=s(x)+s'(P).$

Summarizing (1) and (2), $s(x+P)=s(x)+s'(P)$. \qed
\end{proof}

\begin{remark} Similarly, we can prove that
$s(x-P)=s(x)-s'(P)$, where $s \in \{I,O\},$ $s'=Replace(s),$ $x\in
\mathbb{N}^*,$ $P \in [0]_2$.
\end{remark}

\begin{lemma} \label{lemma:1} (\textbf{Separation Lemma}.)
Suppose $x \in \mathbb{N}^*,$ $s \in \{I,O\}^{\geq 2},$
$s'=Replace(s).$ If $GetS(s',1,j)(P) \in [0]_2$,
$j=0,1,2,...,|s|-1$, then

(1) $IsEven(GetS(s,1,j)(x+P)) = IsEven(GetS(s,1,j)(x))$;

(2) $GetS(s,1,j+1)(x+P)=GetS(s,1,j+1)(x)+GetS(s',1,j+1)(P).$
\end{lemma}

\begin{proof}

(1) $j=0$.

(1.1) $GetS(s',1,j)(P) \in [0]_2$. That is, $GetS(s',1,0)(P)=P \in
[0]_2.$ Thus, $IsEven(x+P)=IsEven(x)$ due to Lemma
\ref{lemma:thefirstct}. Thus, the intermediate next Collatz
transformation of $x+P$ and $x$ are identical.

(1.2) $GetS(s,1,j+1)(x+P) \\
=GetS(s,1,1)(x+P) \; \because j=0\\
=GetS(s,1,1)(x)+GetS(s',1,1)(P) \;\;\;\;\; \because$ Lemma
\ref{lemma:1|STR|=1}\\
$=GetS(s,1,j+1)(x)+GetS(s',1,j+1)(P).$

(2) $j=1$.

(2.1) Due to (1), $GetS(s,1,1)(x+P) =
GetS(s,1,1)(x)+GetS(s',1,1)(P)$.

Besides, $GetS(s',1,1)(P) \in [0]_2$. Thus,

$IsEven(GetS(s,1,1)(x+P)) = IsEven(GetS(s,1,1)(x))$. Thus, the
intermediate next Collatz transformation of $x+P$ and $x$ are
identical.

(2.2) There exists two cases as follows:

(2.2.1) If $GetS(s,1,j+1)=GetS(s,1,j)\|I$, then

$GetS(s,1,j+1)(x+P) \\
=(GetS(s,1,j)\|I)(x+P)\\
=I(GetS(s,1,j)(x+P)) \\
=I(GetS(s,1,j)(x)+GetS(s',1,j)(P)) \;\;\;\;\; \because (1.2) \\
=(3(GetS(s,1,j)(x)+GetS(s',1,j)(P))+1)/2\\
=(3*GetS(s,1,j)(x)+1)/2+3*GetS(s',1,j)(P)/2\\
=I(GetS(s,1,j)(x))+I'(GetS(s',1,j)(P))\\
=(GetS(s,1,j)\|I)(x)+(GetS(s',1,j)\|I')(P)\\
=GetS(s,1,j+1)(x)+GetS(s',1,j+1)(P).$

(2.2.2) If $GetS(s,1,j+1)=GetS(s,1,j)\|O$, then

$GetS(s,1,j+1)(x+P) \\
=(GetS(s,1,j)\|O)(x+P)\\
=O(GetS(s,1,j)(x+P))\\
=O(GetS(s,1,j)(x)+GetS(s',1,j)(P)), \;\;\;\;\; \because (1.2) \\
=(GetS(s,1,j)(x)+GetS(s',1,j)(P))/2\\
=O(GetS(s,1,j)(x))+O(GetS(s',1,j)(P))\\
=(GetS(s,1,j)\|O)(x)+(GetS(s',1,j)\|O)(P)\\
=GetS(s,1,j+1)(x)+GetS(s',1,j+1)(P).$

(Note that, here $j+1=2.$ Recall that ``$\|$'' is concatenation.)

(3) Similarly, $j=2$.

Due to (2), $GetS(s,1,2)(x+P)=GetS(s,1,2)(x)+GetS(s',1,2)(P)$.

Besides, $GetS(s',1,2)(P) \in [0]_2$. Thus,

$IsEven(GetS(s,1,j)(x+P))=IsEven(GetS(s,1,j)(x)).$ Thus, the
intermediate next Collatz transformation of $x+P$ and $x$ are
identical.

(Note that, here $j=2$).

Again, we can prove the following similar to (2.2).

$GetS(s,1,j+1)(x+P)
%=(Substr(s,1,j)\|I)(x+P) \; \vee \; (Substr(s,1,j)\|O)(x+P)\\
%=I(Substr(s,1,j)(x+P)) \; \vee \; O(Substr(s,1,j)(x+P))\\
%=I(Substr(s,1,j)(x)+Substr(s',1,j)(P)) \; \vee \\
%O(Substr(s,1,j)(x)+Substr(s',1,j)(P))\\
%=(3(Substr(s,1,j)(x)+Substr(s',1,j)(P))+1)/2 \; \vee \\
%(Substr(s,1,j)(x)+Substr(s',1,j)(P))/2\\
%=(3*Substr(s,1,j)(x)+1)/2+3*Substr(s',1,j)(P)/2 \; \vee\\
%O(Substr(s,1,j)(x))+O(Substr(s',1,j)(P))\\
%=I(Substr(s,1,j)(x))+I'(Substr(s',1,j)(P)) \; \vee\\
%O(Substr(s,1,j)(x))+O(Substr(s',1,j)(P))\\
%=(Substr(s,1,j)\|I)(x)+(Substr(s',1,j)\|I')(P) \; \vee\\
%(Substr(s,1,j)\|O)(x)+(Substr(s',1,j)\|O)(P)\\
=GetS(s,1,j+1)(x)+GetS(s',1,j+1)(P). \;\;$

(Note that, here $j+1=3.$)

(4) Similarly, we can prove $j=3,...,|s'|-1,$ respectively and
especially in an order. \qed
\end{proof}

\begin{remark}
\item (1) Obviously, above conclusion can be extended to include $|s|=1$ by
Lemma \ref{lemma:thefirstct} and Lemma \ref{lemma:1|STR|=1}.

\item (2) Separation Lemma states the sufficient condition
for guaranteeing that all intermediate parities of transformed
integers $x+P$ are exactly identical with those of $x$.

\item (3) Separation Lemma is general, as $s$ could
be either original dynamics or reduced dynamics. Besides, the length
of $s$ can be omitted by using condition ``$j=0,1,2,...$'' instead
of ``$j=0,1,2,...,|s|-1$''.
\end{remark}

We next explore how to compute $GetS(s',1,j)(x)$.

\begin{definition} Function $CntI(\cdot).$ $CntI: s
\rightarrow n$. It takes as input $s \in \{I,O\}^{\geq1}$, and
outputs $n \in \mathbb{N}$ that is the count of $I$ in $s$.
\end{definition}

%\begin{definition} Function $CntO(\cdot)$. $CntO: s \rightarrow n$
%takes as input $s \in \{I,O\}^{\geq1}$, and outputs $n\in
%\mathbb{N}^*$ that is the count of ``$O$'' in $s$.
%\end{definition}

\begin{example} $CntI(IIOO)=2,$ $CntI(III)=3.$ Obviously, the function name
stems from ``Count the number of $I$''.\end{example}

\begin{lemma} \label{lemma:2}
Suppose $s\in \{I,O\}^{\geq 1},$ $s'=Replace(s),$ $x \in
\mathbb{N}^*$, we have
$$GetS(s',1,j)(x)=\frac{3^{CntI(GetS(s,1,j))}}{2^j}*x, \;\;
j=1,2,...,|s|.$$
\end{lemma}

\begin{proof}
(1) $|s|=1.$ Thus, $j=1.$

(1.1) If $s=I,$ then $s'=Replace(s)=I'.$

$GetS(s',1,j)(x)=GetS(I',1,1)(x)=I'(x)=3*x/2=3^1*x/2^1\\
=3^{CntI(I)}*x/2^{|I|} = 3^{CntI(GetS(s,1,j))}/2^j*x.$

(1.2) If $s=O,$ then $s'=Replace(s)=O.$

$GetS(s',1,j)(x)=O(x)=x/2=3^0*x/2^1=3^{CntI(O)}*x/2^{|O|}\\
= 3^{CntI(GetS(s,1,j))}/2^j*x.$

(2) $|s| \geq 2.$

(2.1) $j=1$.

(2.1.1) If $GetS(s,1,1)=I,$ then

$GetS(s',1,j)(x)=I'(x)=3*x/2=3^1*x/2^1=3^{CntI(I)}*x/2^{|I|}
=3^{CntI(GetS(s,1,j))}/2^j*x.$

(2.1.2) If $GetS(s,1,1)=O,$ then

$GetS(s',1,j)(x)=O(x)=x/2=3^0*x/2^1=3^{CntI(O)}*x/2^{|O|}\\
=3^{CntI(GetS(s,1,j))}/2^j*x.$

(2.2) Iteratively, for $j=1,2,...,|s'|-1$ in an order (recall that
$|s'|=|s|$).

(2.2.1) If $GetS(s,2,1)=I,$ then

$GetS(s',1,j+1)(x) =(GetS(s',1,j)\|I')(x)\\
=I'(GetS(s',1,j)(x)) =I'(GetS(s',1,j)(x))\\
=3*GetS(s',1,j)(x)/2\\
=3*\frac{3^{CntI(GetS(s,1,j))}}{2^j}*x/2 \;\;\; \because (2.1) for \; j=1, (2.2) for \; j=2,...,|s'|-1\\
=\frac{3^{CntI(GetS(s,1,j))+1}}{2^{j+1}}*x \\
=\frac{3^{CntI(GetS(s,1,j+1))}}{2^{j+1}}*x.  \;\;\; \because
GetS(s,1,j+1)=GetS(s,1,j)\|I$

(2.2.2) If $GetS(s,2,1)=O,$ then

$GetS(s',1,j+1)(x) =(GetS(s',1,j)\|O)(x)\\
=O(GetS(s',1,j)(x)) =GetS(s',1,j)(x)/2 \\
=\frac{3^{CntI(GetS(s,1,j))}}{2^j}*x/2 =\frac{3^{CntI(GetS(s,1,j))}}{2^{j+1}}*x \\
=\frac{3^{CntI(GetS(s,1,j+1))}}{2^{j+1}}*x.  \;\;\; \because
GetS(s,1,j+1)=GetS(s,1,j)\|O$ \qed
\end{proof}

\begin{remark} \label{remark:ss'}
Recall that $s'=GetS(s',1,|s|)$ and $s=GetS(s,1,|s|)$. Thus, when
$j=|s|$, then
$s'(x)=GetS(s',1,|s|)(x)=\frac{3^{CntI(GetS(s,1,|s|))}}{2^{|s|}}*x
=\frac{3^{CntI(s)}}{2^{|s|}}*x$.

%\item (2) Indeed, $j=|GetS(s,1,j))|,$ thus the lemma can be restated
%as\\
%$GetS(s',1,j)(x)=\frac{3^{CntI(GetS(s,1,j))}}{2^{|GetS(s,1,j)|}}*x,$
%$j=1,2,...,|s|.$
\end{remark}

Similar to Lemma \ref{lemma:1} (Separation Lemma), a variant of
Separation Lemma can be given as follows:

\begin{lemma} \label{lemma:1extended}
%$s \in \{I,O\}^{\geq 2},$ $s'=Replace(s),$ $x, P, x-P \in
%\mathbb{N}^*$.
%
%(1)
%\[GetS(s,1,1)(x-P)=GetS(s,1,1)(x)-GetS(s',1,1)(P).\]
%
%(2) $j=1,2,...,|s'|-1.$
%\[Substr(STR,1,j)(x-P)=Substr(STR,1,j)(x)-Substr(STR',1,j)(P) \wedge\]
%\[Substr(STR',1,j)(P) \in [0]_2  \Rightarrow\]
%\[IsOdd(Substr(STR,1,j)(x-P)) = IsOdd(Substr(STR,1,j)(x)) \wedge\]
%\[Substr(STR,1,j+1)(x-P)=Substr(STR,1,j+1)(x)-Substr(STR',1,j+1)(P).\]

Suppose $x \in \mathbb{N}^*,$ $s \in \{I,O\}^{\geq 2},$
$s'=Replace(s).$ If $GetS(s',1,j)(P) \in [0]_2$,
$j=0,1,2,...,|s|-1$, then

(1) $IsEven(GetS(s,1,j)(x-P)) = IsEven(GetS(s,1,j)(x))$;

(2) $GetS(s,1,j+1)(x-P)=GetS(s,1,j+1)(x)-GetS(s',1,j+1)(P).$
\end{lemma}
\begin{proof}
The proof is similar to Lemma \ref{lemma:1}. \qed
\end{proof}

\section{Partition Theorem \label{sec:partition}}

\begin{notation} $\mathsf{DYNM}(x,n) \in \{I,O\}^n.$ It
denotes the first $n \in \mathbb{N}^*$ transformations (in terms of
$I$ or $O$) of $x \in \mathbb{N}^*$ no matter whether
$\mathsf{RD}[x]$ exists or not.
\end{notation}

For example, $\mathsf{DYNM}(19,2)=II,$ $\mathsf{DYNM}(19,3)=IIO,$
and $\mathsf{DYNM}(19,4)=IIOO.$

Interestingly, we observe facts as follows:

$\mathsf{DYNM}(x \in [3]_8,3)=IIO$. Let $s=IIO.$

$\forall x \in [3]_8,$ if $s(x) \in [0]_2$, then $s(x+8)\in [1]_2;$
If $s(x)\in [1]_2,$ then $s(x+8)\in [0]_2$.

For example, $x=3 \in [3]_8,$ $s(3)=IIO(3)=IO(5)=O(8)=4 \in [0]_2$,
$s(3+8=11)=IIO(11)=IO(17)=O(26)=13 \in [1]_2.$ Moreover,
$s(11+8=19)=IIO(19)=IO(29)=O(44)=22 \in [0]_2.$

\begin{remark}
For better understanding our further proofs, a graph called reduced
dynamics graph is proposed and explained in Appendix, which can
visualize the following conclusions, but not mandatory for the
proof. The details on tree-based dynamics graph is presented in
another paper \cite{weiispa}.
\end{remark}

Next, we will prove above key observations.

Recall that by \emph{Subset Theorem} (Theorem \ref{th:match}), if
$\mathsf{DYNM}(x \in r_1,n)=s$, $\mathsf{DYNM}(x \in r_2, n+1)=s\|I$
or $s\|O$, then $r_2 \subset r_1$.

Indeed, in this section we will prove that not only $r_2 \subset
r_1,$ but also $r_2$ is either half partition of $r_1$. More
specifically, the first distinction between transformations of $x$
and transformations of $x+2^t$ occurs at the $(t+1)$-th
transformation (if $|\mathsf{RD}[x]|>t+1$). That is, $x \equiv x+2^t
\equiv i \mod 2^t $. $[i]_{2^{t+1}} \cup
[i+2^t]_{2^{t+1}}=[i]_{2^t}$, $[i]_{2^{t+1}} \cap
[i+2^t]_{2^{t+1}}=\emptyset.$ If $x \in [i]_{2^{t+1}}$, then $x+2^t
\in [i+2^t]_{2^{t+1}}.$ If $x \in [i+2^t]_{2^{t+1}}$, then $x+2^t
\in [i]_{2^{t+1}}.$

\begin{lemma} \label{lemma:IsOdd}
If $\mathsf{DYNM}(x,t)=s \in \{I,O\}^t,$ $x,t \in \mathbb{N}^*$,
$j=1,2,...,t-1,$ then \\
(1) $t=1.$ $IsEven(s(x)) \neq IsEven(s(x+2^t))$; Or,\\
(2) $t \geq 2.$
%\begin{equation}
%\left\{ \begin{aligned}
$IsEven(GetS(s,1,j)(x)) = IsEven(GetS(s,1,j)(x+2^t)),$ and
$IsEven(s(x)) \neq IsEven(s(x+2^t))$.
%\end{aligned} \right.
%\end{equation}
\end{lemma}

\begin{proof}

%W.l.o.g., suppose $f^t=I^{p_1}O^{q_1}...I^{p_n}O^{q_n},$
%$\sum_{j=1}^n(p_j+q_j)=t, p_j,q_j \in \mathbb{N}, j=1,2,...,n.$

%Due to Separation Lemma (i.e., Lemma \ref{lemma:1}),

%$f^t(x+2^t)=I^{p_1}O^{q_1}...I^{p_n}O^{q_n}(x+2^t)\\
%=I^{p_1}O^{q_1}...I^{p_n}O^{q_n}(x)+I'^{p_1}O^{q_1}...I'^{p_n}O^{q_n}(2^t)$,
%where $I'^{p_1}O^{q_1}...I'^{p_n}O^{q_n}$ is the replacement of all
%``$I$'' with ``$I'$'' in $I^{p_1}O^{q_1}...I^{p_n}O^{q_n}.$
%
%$I'^{p_1}O^{q_1}...I'^{p_n}O^{q_n}(2^t)\\
%=(\frac{3}{2})^{p_n}...((\frac{3}{2})^{p_2}((\frac{3}{2})^{p_1}(2^t)/2^{q_1})/2^{q_2}).../2^{q_n}\\
%=\frac{3^{\sum_{j=1}^np_j}}{2^{\sum_{j=1}^n(p_j+q_j)}}*2^t=3^p \in
%[1]_2.$
%
%Together with
%$I'^{p_1}O^{q_1}...I'^{p_n}O^{q_n}(2^t)=f^t(x+2^t)-f^t(x)$, thus
%$f^t(x+2^t)-f^t(x) \in [1]_2$. Therefore, $IsOdd(f^t(x)) \neq
%IsOdd(f^t(x+2^t))$.\\

%Besides, let $Substr(f^t,1,k)=f^k, 1 \leq k \leq t-1.$

%$f^k(x+2^t)=I^{p_1}O^{q_1}...I^{p_{k1}}O^{q_{k2}}(x+2^t)\\
%=I^{p_1}O^{q_1}...I^{p_{k1}}O^{q_{k2}}(x)+I'^{p_1}O^{q_1}...I'^{p_{k1}}O^{q_{k2}}(2^t)$,
%where $I'^{p_1}O^{q_1}...I'^{p_{k1}}O^{q_{k2}}$ is the replacement
%of all ``$I$'' with ``$I'$'' in
%$I^{p_1}O^{q_1}...I^{p_{k1}}O^{q_{k2}}.$
%
%$I'^{p_1}O^{q_1}...I'^{p_{k1}}O^{q_{k2}}(2^t)\\
%=(\frac{3}{2})^{p_{k1}}...((\frac{3}{2})^{p_2}((\frac{3}{2})^{p_1}(2^t)/2^{q_1})/2^{q_2}).../2^{q_{k2}}\\
%=\frac{3^{\sum_{j=1}^{k1}p_j}}{2^{\sum_{j=1}^{k1}p_j+\sum_{j=1}^{k2}q_j}}*2^t
%\in [0]_2 \;\;\; \because
%\sum_{j=1}^{k1}p_j+\sum_{j=1}^{k2}q_j<\sum_{j=1}^n(p_j+q_j)=t.$

(1) $t=1.$ $2^t=2 \in [0]_2.$ $s=I$ or $O$.\\
$I(x+2)=(3(x+2)+1)/2=(3x+1)/2+3=I(x)+3.$ \\
Thus, $IsEven(I(x+2)) \neq IsEven(I(x)).$\\
$O(x+2)=(x+2)/2=O(x)+1.$ Thus, $IsEven(O(x+2))\neq IsEven(O(x))$.\\
Hence, $IsEven(s(x)) \neq IsEven(s(x+2^t)).$

(2) $t \geq 2.$ Let $s'=Replace(s)$.

By Lemma \ref{lemma:2},
$GetS(s',1,j)(x)=\frac{3^{CntI(GetS(s,1,j))}}{2^j}*x, \;
j=1,2,...,|s|.$

$GetS(s',1,j)(2^t)=\frac{3^{CntI(GetS(s,1,j))}}{2^j}*2^t
=3^{CntI(GetS(s,1,j))}*2^{t-j} \in [0]_2,$ due to $j \leq t-1$ ($t-j
\geq 1).$

Due to \emph{Separation Lemma} (i.e., Lemma \ref{lemma:1}), we have

(1) $GetS(s',1,j)(x+2^t)=GetS(s,1,j)(x)+GetS(s',1,j)(2^t)$, and

(2) $IsEven(GetS(s,1,j)(x+2^t))=IsEven(GetS(s,1,j)(x)),$ and

(3) $s(x+2^t)=s(x)+s'(2^t).$

Due to Remark \ref{remark:ss'},
$s'(2^t)=\frac{3^{CntI(s)}}{2^{|s|}}*2^t=\frac{3^{CntI(s)}}{2^t}*2^t
=3^{CntI(s')} \in [1]_2.$

Thus, $IsEven(s(x+2^t)) \neq IsEven(s(x)).$ \qed

\end{proof}

\begin{remark} \label{remark:lemmaIsOdd}By assuming $GetS(\cdot,\cdot,0)(x)=x$, above lemma can be restated as
$IsEven(GetS(s,1,j)(x)) = IsEven(GetS(s,1,j)(x+2^t))$ and\\
$IsEven(s(x)) \neq IsEven(s(x+2^t))$ where $j=0,1,...,t-1, t \geq
1.$
\end{remark}

Above lemma states that the first $t$ transformations of $x$ and
$x+2^t$ are identical. It can be extended to include $x-2^t$ upon
$x-2^t>0$. In other words, the first $t$ transformations of $x$ and
$x \pm 2^t$ (i.e., $x \in [x \mod 2^t]_{2^t}$) are identical. Or,
the first $t$ transformations of $x$ is determined by $x \mod 2^t$.
Or, all $x \in [i]_{2^t}$ ($i = x \mod 2^t$) have the same first $t$
transformations.

\begin{lemma} \label{lemma:IsOdd'}
If $\mathsf{DYNM}(x,t)=s \in \{I,O\}^t,$ $x,t \in \mathbb{N}^*$,
then \\
$\mathsf{DYNM}(x \in [i]_{2^t},t)=s,$ where $i = x \mod 2^t$.
\end{lemma}
\begin{proof}

Let $j=0,1,2,...,t-1.$

%Obviously, $IsOdd(x+m*2^t)=IsOdd(x)$, $IsOdd(x-m*2^t)=IsOdd(x),
%x-m*2^t>0.$ Thus, the first transformations for $x$, $x+m*2^t$ and
%$x-m*2^t$ are identical.

(1) By Lemma \ref{lemma:IsOdd} and Remark \ref{remark:lemmaIsOdd},

$IsEven(GetS(s,1,j)(x))=IsEven(GetS(s,1,j)(x+2^t)).$

Thus, $\mathsf{DYNM}(x+2^t,t)=\mathsf{DYNM}(x,t)=s.$

Similarly,

$IsEven(GetS(s,1,j)(x+2^t))=IsEven(GetS(s,1,j)(x+2^t+2^t)).$ Thus,

$\mathsf{DYNM}(x+2^t+2^t,t)=\mathsf{DYNM}(x+2^t,t)=s.$

Iteratively and hence,

$\mathsf{DYNM}(x+m*2^t,t)=\mathsf{DYNM}(x,t)=s, \;\; m \in
\mathbb{N}^*.$

(2) Due to the variant of Separation Lemma (i.e., Lemma
\ref{lemma:1extended}),

$IsEven(GetS(s,1,j)(x))=IsEven(GetS(s,1,j)(x-2^t)).$

When $x-2^t>0$, we have

$\mathsf{DYNM}(x-2^t,t)=\mathsf{DYNM}(x,t)=s.$

Iteratively, when $x-m*2^t>0, m \in \mathbb{N}^*$, we have

$\mathsf{DYNM}(x-m*2^t,t)=\mathsf{DYNM}(x,t)=s.$

Therefore, due to (1) and (2),

$\mathsf{DYNM}(x \in [i]_{2^t},t)=s,$ where $i=x \mod 2^t$. \qed
\end{proof}

%Indeed, above Lemma \ref{lemma:IsOdd0} and Lemma \ref{lemma:IsOdd}
%are also maintained for $i \in [0]_2.$ However, we mainly are
%interested in $i \in [1]_2$ (i.e., $x \in [1]_2$).

\begin{remark}
\item (1) Above lemma includes the special case $i=0$ or $x \in [0]_2$.
However, $\mathsf{DYNM}(x \in [0]_2,1)=\mathsf{RD}[x\in[0]_2]=O$ is
trivial, $\mathsf{DYNM}(x \in [1]_2,n \in \mathbb{N}^*)$ is thus
of more interest. That is, we mainly concentrate on \\
$\mathsf{DYNM}(x \in [i]_{2^t},n \in \mathbb{N}^*), i,t \in
\mathbb{N}^*, 1 \leq i \leq 2^t-1, i \in [1]_2.$

\item (2) In above lemmas (Lemma
\ref{lemma:IsOdd} and Lemma \ref{lemma:IsOdd'}), $s$ is general.
That is, both reduced dynamics and original dynamics satisfy above
lemmas.

\item (3) Recall that, Lemma \ref{lemma:IsOdd'} states that the first
$t$ transformations for $x \in [x \mod 2^t]_{2^t}$ are identical.
Lemma \ref{lemma:IsOdd} states that the $(t+1)$-th transformation
for $x$ and $x+2^t$ is distinct (so-called ``forking''). Note that,
it can be observed that either $x$ or $x+2^t$ falls in either
partition of $[x \mod 2^t]_{2^t}$ (i.e., $[x \mod 2^t]_{2^{t+1}},
[(x \mod 2^t)+2^t]_{2^{t+1}}$), respectively. Moreover, all natural
numbers in $x \in [x \mod 2^t]_{2^t}$ are further partitioned into
two halves, and either partition results in either transformation
(i.e., ``$I$'' or ``$O$'') due to ``forking'', iteratively.

\end{remark}

%\begin{lemma} \label{lemma:powerof2}
%The level of partitions (in terms of residue modules) is the power
%of 2. In other words, the partitions always have the form like $x
%\in [0]_2 \cup [i]_{2^t}, t\in \mathbb{N}^*, 1 \leq i \leq 2^t-1, i
%\in [1]_2, i = x \mod 2^t.$
%\end{lemma}
%
%\begin{proof} Straightforward. As each time one partition is divided into two halves
%(or each residue class is partitioned into two halves), the level of
%partitions in terms of residue modules is always the power of 2.
%
%For $x\in[0]_2,$ the further partition is not of interest because
%$\mathsf{CODE}(x\in [0]_2)=O$ is final.
%
%For $x \in [i]_{2^t},$ the further partitions are $x \in
%[i]_{2^{t+1}}$ and $x\in [i+2^t]_{2^{t+1}}$, as either partition
%results in either transformation in the next (recall Lemma
%\ref{lemma:IsOdd}). More specifically, to determine whether $f^t(x)
%\in [1]_2$ or not, $x \in [i]_{2^t}$ is partitioned into two halves
%- $x\in [i]_{2^{t+1}}$ and $x \in [i+2^t]_{2^{t+1}}$.
%$[i]_{2^{t+1}}\cap [i+2^t]_{2^{t+1}}=\emptyset,$ $[i]_{2^{t+1}} \cup
%[i+2^t]_{2^{t+1}}=[i]_{2^t}$. Therefore, they always have the form
%like $[i]_{2^t},$ $t\in \mathbb{N}^*,$ $1 \leq i \leq 2^t-1,$ $i \in
%[1]_2.$ \qed
%\end{proof}

Similar to Lemma \ref{lemma:IsOdd} and Lemma \ref{lemma:IsOdd'} but
more specifically for whether ``forking'' transformation is ``$I$''
or ``$O$'' and why, we have following theorem. Roughly speaking,
partition residue class $[i]_{2^t}$ determines the first $t$
transformations. If current integer is less than the starting
integer, then reduced dynamics is obtained. %(Indeed, we proved that
%$CntO(s)=\lceil \log_21.5*CntI(s) \rceil, t \geq 2$ or $CntO(s)=1,
%t=1$).
Otherwise, further transformation occurs, and whether it is
``$I$'' or ``$O$'' is determined by and only by either further half
partition of the current residue class of $x$.

\begin{theorem} \label{th:partition}
(\textbf{Partition Theorem}.)\\
(1) $\mathsf{DYNM}(x \in [0]_2,1)=\mathsf{RD}[x \in [0]_2]=O.$\\
(2) Suppose $\mathsf{DYNM}(i,t)=s \in \{I,O\}^t,$ $t,i \in
\mathbb{N}^*,$ $1 \leq i \leq 2^t-1$. We have \\
(2.1) $\mathsf{DYNM}(x \in [i]_{2^t},t)=s.$\\
(2.2) If %$CntO(f^t) = \lceil \lambda*CntI(f^t) \rceil, CntO(s) <
%\lceil \lambda*CntI(s) \rceil,$ $s=Substr(f^t,1,k),$
%$k=1,2,...,t-1,$ $t \geq 2,$
$s(x)<x,$
then $\mathsf{RD}[x \in [i]_{2^t}]=s$.\\
(2.3) If %$CntO(f^t)<\lceil \lambda*CntI(f^t) \rceil, t\geq 2$, then \\
$s(x)\not <x$ and $s(i) \in [0]_2$, then\\
$\mathsf{DYNM}(x \in [i]_{2^{t+1}},t+1)=s\|O$ and
$\mathsf{DYNM}(x \in [i+2^t]_{2^{t+1}},t+1)=s\|I$;\\
If $s(x)\not <x$ and $s(i) \in [1]_2$, then\\
$\mathsf{DYNM}(x \in [i]_{2^{t+1}},t+1)=s\|I$ and $\mathsf{DYNM}(x
\in [i+2^t]_{2^{t+1}},t+1)=s\|O$.
\end{theorem}

\begin{proof}
(1) Straightforward.

(2.1) The proof is similar to Lemma \ref{lemma:IsOdd} (2) (except
that here $t=1$ is possible), or by Lemma \ref{lemma:IsOdd'} (1).

%Suppose $\mathsf{DYNM}(i,t)=f^t.$ $1 \leq i \leq 2^t-1, t \in
%\mathbb{N}^*$. Next, we will prove
%$\mathsf{DYNM}(i+m*2^t,t)=f^t$ where $m \in \mathbb{N}^*.$
%
%Let $f'^k = Replace(f^k),$ $f^k=Substr(f^t,1,k),k=0,1,...,t-1.$
%Recall that, $f^0(x)=Substr(\cdot,\cdot,0)(x)=f'^0(x)=x.$
%
%For $k=0$, $f^k(i+m*2^t)=i+m*2^t=f^k(i)+f'^k(m*2^t).$ \\
%Besides, $m*2^t \in [0]_2 \Rightarrow IsOdd(i+m*2^t)=IsOdd(i) \\
%\Rightarrow IsOdd(f^k(i+m*2^t))=IsOdd(f^k(i)), k=0.$
%
%For $k \geq 1,$ due to Lemma \ref{lemma:2'},\\
%$f'^k(m*2^t)=\frac{3^{CntI'(f'^k)}}{2^{|f'^k|}}*m*2^t=\frac{3^{CntI'(f'^k)}}{2^k}*m*2^t\\
%=3^{CntI'(f'^k)}*m*2^{t-k}  \in [0]_2, \;\; \because m \in
%\mathbb{N}^*, 1 \leq k \leq t-1 \Rightarrow t-k \geq 1.$
%
%Due to Separation Lemma (i.e., Lemma \ref{lemma:1}), we have
%
%$f^k(i+m*2^t)=f^k(i)+f'^k(m*2^t), k=0,1,...,t-1 \;\; \wedge$
%
%$IsOdd(f^k(i+m*2^t))=IsOdd(f^k(i)), k=0,1,...,t-1 \;\; \wedge$
%
%$f^t(i+m*2^t)=f^t(i)+f'^t(m*2^t).$
%
%Recall that $k=0$ means the first transformation (i.e., $I$ or $O$)
%of $i+m*2^t$ and $i$ are identical.
%
%Therefore, $\mathsf{DYNM}(x \in
%[i]_{2^t},t)=\mathsf{DYNM}(i,t)=f^t$.

(2.2) %By Corollary \ref{th:codetheoremsufficiency}, if
%$CntO(f^t)=\lceil \lambda*CntI(f^t) \rceil, t\geq 2$,
If $s(x)<x,$ then reduced dynamics ends and $\mathsf{RD}[x \in
[i]_{2^t}]=\mathsf{DYNM}(x \in [i]_{2^t},t)=s$.

(2.3) We inspect the value of $IsEven(s'(m*2^t))$ where
$s'=Replace(s)$ to decide whether $IsEven(s(i+m*2^t)) = IsEven(s(i))$ or not,\\
since $s(i+m*2^t)=s(i)+s'(m*2^t)$, by Separation Lemma (i.e.,
Lemma \ref{lemma:1}).

By Remark \ref{remark:ss'} (recall that $t \in \mathbb{N}^*$),

$s'(m*2^t)=\frac{3^{CntI(s)}}{2^{|s|}}*m*2^t=\frac{3^{CntI(s)}}{2^t}*m*2^t
=3^{CntI(s)}*m.$

$m \in [1]_2 \Rightarrow 3^{CntI(s)}*m \in [1]_2  \Rightarrow
IsEven(s(i+m*2^t)) \neq IsEven(s(i)).$

Otherwise, $m \in [0]_2 \Rightarrow IsEven(s(i+m*2^t)) =
IsEven(s(i)).$

(1) $m \in [1]_2 \Leftrightarrow
x=i+(m-1)*2^t+2^t=i+\frac{m-1}{2}*2^{t+1}+2^t \in
[i+2^t]_{2^{t+1}}.$

Thus, $x \in [i+2^t]_{2^{t+1}} \Rightarrow IsEven(s(x)) \neq
IsEven(s(i)).$

(2) $m \in [0]_2  \Leftrightarrow x=i+m*2^t=i+\frac{m}{2}*2^{t+1}
\in [i]_{2^{t+1}}\setminus\{i\}.$ \\
(Recall that $[0]_2=\{a|a\in \mathbb{N}^*, a \equiv 0 \mod 2\},$
thus $m\geq 2.$)

Thus, $x \in [i]_{2^{t+1}}\setminus\{i\} \Rightarrow
IsEven(s(x))=IsEven(s(i)).$

Obviously, $x=i$ can be included too, because \\
$x=i \Rightarrow IsEven(s(x))=IsEven(s(i))$. Therefore, \\
$x \in [i]_{2^{t+1}} \Rightarrow IsEven(s(x))=IsEven(s(i)).$

(3) Finally, if $s(i) \in [0]_2$, then $\mathsf{DYNM}(x \in
[i]_{2^{t+1}},t+1)=s\|O$ due to (2), and $\mathsf{DYNM}(x \in
[i+2^t]_{2^{t+1}},t+1)=s\|I$ due to (1);

If $s(i) \in [1]_2$, then $\mathsf{DYNM}(x \in
[i]_{2^{t+1}},t+1)=s\|I$ due to (2), and $\mathsf{DYNM}(x \in
[i+2^t]_{2^{t+1}},t+1)=s\|O$ due to (1).

Note that, $[i+2^t]_{2^{t+1}} \cap [i]_{2^{t+1}} = \emptyset,$ and
$[i+2^t]_{2^{t+1}} \cup [i]_{2^{t+1}} = [i]_{2^t}.$
%
%Last, it is worth to stress that above proof is valid (iteratively)
%for all $t \in \mathbb{N}^*$.
\qed
\end{proof}

\begin{remark} \label{remark:partition}

\item (1) We call the theorem as ``Partition Theorem'' because
$[i]_{2^{t}}$ is partitioned into two halves for deciding the
$(t+1)$-th transformation. Recall that, $[i]_{2^t}=[i]_{2^{t+1}}
\cup [i+2^t]_{2^{t+1}}, [i]_{2^{t+1}} \cap
[i+2^t]_{2^{t+1}}=\emptyset.$

\item (2) The level of partitions (in terms of residue modules) is the
power of 2. In other words, the partitions always have the form like
$x \in [0]_2 \cup [i]_{2^t}, t\in \mathbb{N}^*, 1 \leq i \leq 2^t-1,
i \in [1]_2, i = x \mod 2^t.$ \\
For $x \in [0]_2,$ the further partition is not of interest because
$\mathsf{RD}[x\in [0]_2]=O$ is obtained. \\
For $x \in [i]_{2^t},$ the further partitions are $x \in
[i]_{2^{t+1}}$ and $x\in [i+2^t]_{2^{t+1}}$, as either partition
results in either transformation in the next (recall Lemma
\ref{lemma:IsOdd}). More specifically, to determine whether $s(x)
\in [1]_2$ or not, $x \in [i]_{2^t}$ is partitioned into two halves
- $x\in [i]_{2^{t+1}}$ and $x \in [i+2^t]_{2^{t+1}}$.
$[i]_{2^{t+1}}\cap [i+2^t]_{2^{t+1}}=\emptyset,$ $[i]_{2^{t+1}} \cup
[i+2^t]_{2^{t+1}}=[i]_{2^t}$. Therefore, they always have the form
like $[i]_{2^t},$ $t\in \mathbb{N}^*,$ $1 \leq i \leq 2^t-1,$ and
especially, $i \in [1]_2.$

\item (3) This theorem also reveals the link between $\mathsf{RD}[x]$
and $\mathsf{DYNM}(x,t)$ as follows: Suppose $\mathsf{DYNM}(x,t)=s
\in \{I,O\}^t$. Once %$CntO(f^t)=\lceil \lambda*CntI(f^t)\rceil,$
%$CntO(Substr(f^t,1,k))=\lceil \lambda*CntI(Substr(f^t,1,k))\rceil,
%k=1,2,...,t-1,$ $t \geq 2$ or $CntI(f^t)=1,t=1,$
$s(x)<x$, then $\mathsf{DYNM}(x,t)=\mathsf{RD}[x]=s$ (which means
current transformed number has already been less than starting
number for the \emph{first} time); Otherwise, %$CntO(f^t)<\lceil
%\lambda*CntI(f^t)\rceil,$ then
$\mathsf{DYNM}(x,t+1)=s\|I$ or $\mathsf{DYNM}(x,t+1)=s\|O$.

\end{remark}

Next corollary
states that $i \in [0,2^t-1]$ determines the first $t$
transformations of all $x \in [i]_{2^t}.$

\begin{corollary} \label{cor:rstody1}
Suppose $t,i_1,i_2 \in \mathbb{N}^*,$ $0 \leq
i_1,i_2 \leq 2^t-1$. \\
$i_1 = i_2$, if and only if $\mathsf{DYNM}(x \in [i_1]_{2^t},t) =
\mathsf{DYNM}(x \in [i_2]_{2^t},t).$
\end{corollary}
\begin{proof}
(1) Necessity.

(1.1) $t \geq 2,$ thus $1 \leq i_1,i_2 \leq 2^t-1$ due to Theorem
\ref{th:partition} (2) and Remark \ref{remark:partition} (2).

$i_1=i_2 \Rightarrow \mathsf{DYNM}(i_1,t)=\mathsf{DYNM}(i_2,t)\\
\Rightarrow \mathsf{DYNM}(x \in [i_1]_{2^t},t) = \mathsf{DYNM}(x \in
[i_2]_{2^t},t)\\
\because \mathsf{DYNM}(i_1,t)=\mathsf{DYNM}(x \in
[i_1]_{2^t},t), \mathsf{DYNM}(i_2,t) = \mathsf{DYNM}(x \in
[i_1]_{2^t},t)$

(1.2) $t=1,$ thus $0 \leq i_1,i_2 \leq 2^t-1=1.$

$i_1=i_2=(0 \vee 1) \\
\Rightarrow \mathsf{DYNM}(x \in [i_1]_{2^t},t) = \mathsf{DYNM}(x \in [i_1]_2,1)\\
=\mathsf{DYNM}(x \in [i_2]_2,1)=\mathsf{DYNM}(x \in [i_2]_{2^t},t).$

(2) Sufficiency. It is proved by proving converse negative
proposition.

(2.1) $t\geq 2$. Thus $1 \leq i_1,i_2 \leq 2^t-1$.

$i_1 \neq i_2,$ $1 \leq i_1,i_2 \leq 2^t-1$ $\Rightarrow \mathsf{DYNM}(i_1,t) \neq \mathsf{DYNM}(i_2,t) \\
\Rightarrow \mathsf{DYNM}(x \in [i_1]_{2^t},t) \neq \mathsf{DYNM}(x \in[i_2]_{2^t},t). \\
\because \mathsf{DYNM}(x \in [i_1]_{2^t},t)=\mathsf{DYNM}(i_1,t),
\mathsf{DYNM}(x \in [i_2]_{2^t},t)=\mathsf{DYNM}(i_2,t).$

(2.2) $t=1.$

$i_1 \neq i_2 \Rightarrow (i_1=0,i_2=1) \vee (i_1=1,i_2=0)\\
\Rightarrow \mathsf{DYNM}(x \in [i_1]_{2^t},t) = \mathsf{DYNM}(x \in [i_1]_2,1) \\
\neq \mathsf{DYNM}(x \in [i_2]_2,1)=\mathsf{DYNM}(x \in
[i_2]_{2^t},t).$  \qed
\end{proof}

In above corollary, $t$ is categorized for tackling the case $i=0$
to avoid $\mathsf{DYNM}(i=0,t)$.

%Next corollary can be obtained directly by above corollary, but we
%still give a formal proof independent to above corollary and only by
%Partition Theorem.

\begin{corollary} \label{cor:rstody2}
Suppose $x_1,x_2,t \in \mathbb{N}^*, x_1 \neq x_2$.\\
$\mathsf{DYNM}(x_1,t)=\mathsf{DYNM}(x_2,t)$, if and only if $x_1
\equiv x_2 \mod 2^t.$
\end{corollary}
\begin{proof} Straightforward due to Corollary \ref{cor:rstody1}.
\qed
%(1) ($\Rightarrow$). It is proved by proving converse negative
%proposition.
%
%(1.1) $t\geq 2.$ $x_1 \not \equiv x_2 \mod 2^t \\
%\Rightarrow x_1 \in [i_1]_{2^t}, x_2 \in [i_2]_{2^t}, i_1 \neq i_2, 1 \leq i_1,i_2 \leq 2^t-1 \\
%\Rightarrow \mathsf{DYNM}(x_1 \in [i_1]_{2^t},t)= \mathsf{DYNM}(i_1,t) \\
%\neq \mathsf{DYNM}(i_2,t) = \mathsf{DYNM}(x_2 \in [i_2]_{2^t},t)\\
%\Rightarrow \mathsf{DYNM}(x_1,t) \neq \mathsf{DYNM}(x_2,t).$
%
%(1.2) $t=1.$ $\mathsf{DYNM}(x_1,t)=\mathsf{DYNM}(x_2,t)=(I \vee O), x_1 \neq x_2 \\
%\Rightarrow (x_1,x_2 \in [1]_2) \vee (x_1,x_2 \in [0]_2) \Rightarrow
%x_1 \equiv x_2 \mod 2^1.$
%
%(2) ($\Leftarrow$).
%
%(2.1) $t\geq 2.$ (Thus $x_1,x_2 \in [1]_2.$)
%
%$x_1 \equiv x_2 \mod 2^t \Rightarrow x_1,x_2 \in [i]_{2^t},
%i=x_1\mod 2^t=x_2\mod 2^t \\
%\Rightarrow \mathsf{DYNM}(x_1,t) = \mathsf{DYNM}(i,t) =
%\mathsf{DYNM}(x_2,t).$
%
%(2.2) $t = 1.$ $x_1 \equiv x_2 \mod 2^1 \Rightarrow x_1,x_2 \in
%[1]_2 \vee x_1,x_2 \in [0]_2 \\
%\Rightarrow \mathsf{DYNM}(x_1,t)=\mathsf{DYNM}(x_2,t)=(I
%\vee O).$ \qed
\end{proof}

%Corollary \ref{cor:rstody1} and Corollary \ref{cor:rstody2} both
%imply that there exists bijective mapping between paths with length
%$t$ in code graph (i.e., $f^t_i$) and residue classes with module
%$2^t$ (i.e., $[i]_{2^t}$). Corollary \ref{cor:rstody1} emphasizes
%that residue classes determine paths with length $t$ (and more sense
%on residue $i$). Corollary \ref{cor:rstody2} stresses that paths
%with length $t$ determine residue classes (and more sense on
%individual $x$).

Next corollary states that if the first different transformation of
$x_1,x_2$ occurs at the $(t+1)$-th transformation, $t \in
\mathbb{N}^*$, then $x_1 \equiv x_2 \mod 2^t$ and $x_1 \equiv x_2 +
2^t \mod 2^{t+1}.$ As $\mathsf{DYNM}(x\in [0]_2,1)=\mathsf{RD}[x\in
[0]_2]=O$, only $x\in [1]_2$ is of interests.

\begin{corollary} \label{cor:forking} (\textbf{Forking Corollary}.)
Suppose $x_1,x_2 \in [1]_2, x_1\neq x_2, t \in \mathbb{N}^*$.\\
$\mathsf{DYNM}(x_1,t) = \mathsf{DYNM}(x_2,t)$ and
$\mathsf{DYNM}(x_1,t+1) \neq \mathsf{DYNM}(x_2,t+1)$, if and only if
$x_1 \equiv x_2 \mod 2^t$ and $x_1 \equiv x_2+2^t \mod 2^{t+1}.$
\end{corollary}
\begin{proof}
It is straightforward due to Corollary \ref{cor:rstody2}. %We still
%give a formal proof independent to Corollary \ref{cor:rstody2} as
%follows:
%
%(1) ($\Rightarrow$).
%
%(1.1) $x_1 \not \equiv x_2 \mod 2^t \\
%\Rightarrow \mathsf{DYNM}(x_1,t)= \mathsf{DYNM}(x_1 \mod
%2^t,t) \;\;\; \because$ Lemma \ref{lemma:IsOdd'} \\
%$\neq \mathsf{DYNM}(x_2 \mod 2^t,t) =\mathsf{DYNM}(x_2,t). $
%
%Thus, $\mathsf{DYNM}(x_1,t) = \mathsf{DYNM}(x_2,t) \Rightarrow x_1
%\equiv x_2 \mod 2^t.$
%
%(1.2) $x_1 \equiv x_2 \mod 2^{t+1}\\
%\Rightarrow \mathsf{DYNM}(x_1,t+1) = \mathsf{DYNM}(x_1 \mod 2^{t+1},t+1) \\
%= \mathsf{DYNM}(x_2 \mod 2^{t+1},t+1)=\mathsf{DYNM}(x_2,t+1).$
%
%Thus,
%
%$\mathsf{DYNM}(x_1,t+1) \neq \mathsf{DYNM}(x_2,t+1) \\
%\Rightarrow x_1 \not \equiv x_2 \mod 2^{t+1} \\
%\Rightarrow x_1 \equiv x_2+2^t \mod 2^{t+1}. \;\;\; \because (1) \;
%x_1 \equiv x_2 \mod 2^t$
%
%(2) ($\Leftarrow$).
%
%$x_1 \equiv x_2 \mod 2^t \Rightarrow x_1, x_2 \in [i]_{2^t}, i = x_1 \mod 2^t = x_2 \mod 2^t \\
%\Rightarrow \mathsf{DYNM}(x_1,t) = \mathsf{DYNM}(i,t) =
%\mathsf{DYNM}(x_2,t). \;\; \because $ Theorem \ref{th:partition} (2)
%
%$x_1 \equiv x_2+2^t \mod 2^{t+1} \Rightarrow x_1 \not \equiv x_2 \mod 2^{t+1} \\
%\Rightarrow \mathsf{DYNM}(x_1,t+1) = \mathsf{DYNM}([x_1 \mod 2^{t+1}]_{2^{t+1}},t+1) \\
%\neq \mathsf{DYNM}([x_2 \mod 2^{t+1}]_{2^{t+1}},t+1) =
%\mathsf{DYNM}(x_2,t+1). \;\; \because$ Theorem \ref{th:partition}
%(2)
\qed
\end{proof}

\begin{corollary} \label{cor:typereason} (\textbf{Extended Forking Corollary}.)
Suppose \\
$\mathsf{DYNM}(x \in [i_p]_{2^t},t)=s_p \in \{I, O\}^t$,
$p=1,...,m$. \\
$\mathsf{DYNM}(x \in [j_q]_{2^{t+1}},t+1)=s_q \in \{I,O\}^{t+1}$,
$q=1,2,...,n.$ \\
$t\geq 3, t \in \mathbb{N}^*.$
$\forall p \in [1,m],$ $\exists q \in [1,n]$ such that\\
$(s_q=s_p\|I, s_{q+1}=s_p\|O)$ $\vee$ $(s_q=s_p\|O,
s_{q+1}=s_p\|I).$\\
We have following conclusions:\\\\
(1) $\bigcup_{q=1}^n[j_q]_{2^{t+1}}=\bigcup_{p=1}^m[i_p]_{2^t}.$
$i_p \in [1]_2.$ $j_q \in [1]_2.$ $n=2*m.$\\
$s_p=GetS(s_q,1,t)=GetS(s_{q+1},1,t),$ where $q=2*p-1$.\\
(2) $j_{2*p-1}=i_p, j_{2*p}=i_p+2^t.$\\%, p=1,2,...,m.$ Or,
%$j_q=i_{(q+1)/2},$ $j_{q+1}=j_q+2^t, q=1,3,...,2*m-1.$\\
(3) $IsEven(s_p(x \in [i_p]_{2^{t+1}})) \neq IsEven(s_p(x \in
[i_p+2^t]_{2^{t+1}})).$\\
(4) $\mathsf{DYNM}(x \in [i_p]_{2^{t+1}}, t+1) =s_p\|I
\Leftrightarrow \mathsf{DYNM}(x \in
[i_p+2^t]_{2^{t+1}},t+1)=s_p\|O.$\\
(5) $\mathsf{DYNM}(x \in [i_p]_{2^{t+1}}, t+1)=s_p\|O
\Leftrightarrow \mathsf{DYNM}(x \in [i_p+2^t]_{2^{t+1}},
t+1)=s_p\|I.$\\
%(3') $IsOdd(f^t_{(q+1)/2}(x \in [j_q]_{2^{t+1}})) \neq
%IsOdd(f^t_{(q+1)/2}(x \in [j_{q+1}]_{2^{t+1}})), q=1,3,...,2*m-1.$
%
%Or, $\mathsf{DYNM}(x \in [j_q]_{2^{t+1}}, t+1)=f^t_{(q+1)/2}\|I, q=1,3,...,2*m-1\\
%\Leftrightarrow \mathsf{DYNM}(x \in
%[j_{q+1}]_{2^{t+1}},t+1)=f^t_{(q+1)/2}\|O.$
%
%Or, $\mathsf{DYNM}(x \in [j_q]_{2^{t+1}}, t+1)=f^t_{(q+1)/2}\|O, q=1,3,...,2*m-1\\
%\Leftrightarrow \mathsf{DYNM}(x \in
%[j_{q+1}]_{2^{t+1}},t+1)=f^t_{(q+1)/2}\|I.$
(6) $s_p(x \in [i_p]_{2^t}) \in [1]_2 \Leftrightarrow
\mathsf{DYNM}(x \in [i_p]_{2^{t+1}},t+1)=s_p\|I.$\\
$s_p(x \in [i_p]_{2^t}) \in [0]_2\Leftrightarrow
\mathsf{DYNM}(x \in [i_p]_{2^{t+1}},t+1)=s_p\|O.$\\
%(4')
%
%$f^t_p(x \in [i_p]_{2^t}) \in [1]_2, p \in [1,m] \\
%\Leftrightarrow \mathsf{DYNM}(x \in
%[j_{2*p-1}]_{2^{t+1}},t+1)=f^t_p\|I.$
%
%$f^t_p(x \in [i_p]_{2^t}) \in [0]_2, p \in [1,m] \\
%\Leftrightarrow \mathsf{DYNM}(x \in
%[j_{2*p-1}]_{2^{t+1}},t+1)=f^t_p\|O.$
%
%(4'')
%
%$f^t_p(x \in [i_p]_{2^t}) \in [1]_2, p \in [1,m] \\
%\Leftrightarrow \mathsf{DYNM}(x \in
%[j_{2*p}]_{2^{t+1}},t+1)=f^t_p\|O.$
%
%$f^t_p(x \in [i_p]_{2^t}) \in [0]_2, p \in [1,m] \\
%\Leftrightarrow \mathsf{DYNM}(x \in
%[j_{2*p}]_{2^{t+1}},t+1)=f^t_p\|I.$
(7) $\|\{q|GetS(s_q,1,t)(x \in [j_q]_{2^{t+1}}) \in [1]_2\}\|=n/2=m;\\
\|\{q|GetS(s_q,1,t)(x \in [j_q]_{2^{t+1}}) \in [0]_2\}\|=n/2=m.$\\
$GetS(s_q,1,t)(x \in [j_q]_{2^{t+1}}) \in [1]_2
\Leftrightarrow GetS(s_q,1,t)(x \in [j_{q+1}]_{2^{t+1}}) \in [0]_2.$\\
$GetS(s_q,1,t)(x \in [j_q]_{2^{t+1}}) \in [0]_2 \Leftrightarrow
GetS(s_q,1,t)(x \in [j_{q+1}]_{2^{t+1}}) \in [1]_2.$

\end{corollary}

\begin{proof}

By Theorem \ref{th:partition} (Partition Theorem),

if $s_p(x \in [i_p]_{2^t}) \in [1]_2,$ $p \in [1,m]$, then $\exists
q=2*p-1 \in [1,n]$, such that $\mathsf{DYNM}(x \in
[j_q]_{2^{t+1}},t+1)=s_p\|I=s_q.$ Besides, $j_q=i_p$ due to
$IsMatched(GetS(s_q,1,t)(i_p) \in [1]_2, GetS(s_q,t+1,1)=I)=True$.

If $s_p(x \in [i_p]_{2^t}) \in [0]_2,$ $p \in [1,m]$, then $\exists
q=2*p-1 \in [1,n]$, such that $\mathsf{DYNM}(x \in
[j_q]_{2^{t+1}},t+1)=s_p\|O=s_q.$ Besides, $j_q=i_p$ due to
$IsMatched(GetS(s_q,1,t)(i_p) \in [0]_2, GetS(s_q,t+1,1) = O)=True$.

If $s_p(x \in [i_p]_{2^t}) \in [0]_2,$ $p \in [1,m]$, then $\exists
q=2*p \in [1,n]$, such that $\mathsf{DYNM}(x \in
[j_q]_{2^{t+1}},t+1)=s_p\|I=s_q.$ Besides, $j_q=i_p+2^t$ due to
$IsMatched(GetS(s_q,1,t)(j_q=i_p+2^t) \in [1]_2,
GetS(s_q,t+1,1)=I)=True$.

If $s_p(x \in [i_p]_{2^t}) \in [1]_2,$ $p \in[1,m]$, then $\exists
q=2*p \in [1,n]$, such that $\mathsf{DYNM}(x \in
[j_q]_{2^{t+1}},t+1)=s_p\|O=s_q.$ Besides, $j_q=i_p+2^t$ due to
$IsMatched(GetS(s_q,1,t)(j_q=i_p+2^t) \in [0]_2,
GetS(s_q,t+1,1)=O)=True$.

Obviously, $s_p=GetS(s_q,1,t),$ as $s_p\|\{I,O\}=s_q$ where
$q=2*p-1$ or $2*p$ for $p=1,2,...,m.$

$t\geq 3$, thus $i_p \in [1]_2, p=1,2,...,m$. $j_q=i_p \in [1]_2$
where $q=2*p-1,$ and $j_q=i_p+2^t \in [1]_2$ where $q=2*p.$
%Obviously, $p=(q+1)/2$ where $q=1,3,...,2*m-1.$

$p=1,2,...,m$, thus $q=1,2,...,n$ and $n=2*m.$

Besides, $\forall p \in [1,m], [i_p]_{2^t} = [i_p]_{2^{t+1}} \cup
[i_p+2^t]_{2^{t+1}},$ $[i_p]_{2^{t+1}}=[j_q]_{2^{t+1}}$ where
$q=2*p-1,$ $[i_p+2^t]_{2^{t+1}}=[j_q]_{2^{t+1}}$ where $q=2*p.$
Therefore,
$\bigcup_{q=1}^n[j_q]_{2^{t+1}}=\bigcup_{p=1}^m[i_p]_{2^t}.$

In summary, (1) is proved.

Besides, $j_{2*p-1}=i_p,$ $j_{2*p}=i_p+2^t$ where $p=1,2,...,m.$\\
$2*p-1=q$, thus $j_q=i_p=i_{(q+1)/2}$ where $q=1,3,...,2*m-1.$\\
$j_{q+1}=j_{2*p}=i_p+2^t=j_q+2^t$ where $q=1,3,...,2*m-1.$

Thus, (2) is proved.

(3),(4),(5) are due to Lemma \ref{lemma:IsOdd}. %Lemma \ref{lemma:IsOdd0} and

%(3') is the same as (3), as $j_q=i_p, j_{q+1}=i_p+2^t,
%q=1,3,...,2*m-1, p=(q+1)/2$ due to (2).

(6) is due to Partition Theorem, %Remark \ref{remark:powerof2} (2),
and the definition of $\mathsf{DYNM}(\cdot,\cdot)$.
%
%(4') is the same as (4) due to $j_{2*p-1}=i_p$ in (2).
%
%(4'') is the same as (4) due to $j_{2*p}=i_p+2^t$ in (2) and (3) (or
%(3')).

(7) is due to (2).% and (3').
\qed
\end{proof}

Roughly speaking, above corollary states the general (iteratively
and accumulatively) effect of Forking Corollary (Corollary
\ref{cor:forking}). %It can be looked as Extended Forking Corollary.

\begin{lemma} Given $\forall n \in \mathbb{N}^*$, $\exists x \in \mathbb{N}^*$
such that $\mathsf{DYNM}(x,n)=I^n$.
\end{lemma}
\begin{proof}
Straightforward. Given $n \in \mathbb{N}^*$, $\exists x=2^n-1 \in
\mathbb{N}^*$, such that $\mathsf{DYNM}(x,n)=I^n$. Indeed, given
$\forall n \in \mathbb{N}^*$, $\exists x \in [2^n-1]_{2^n}$ such
that $\mathsf{DYNM}(x,n)=I^n$. \qed
\end{proof}

\begin{lemma}
$\forall n \in \mathbb{N}^*,$ $\exists x=2^n-1 \in \mathbb{N}^*$
such that $\mathsf{DYNM}(x,n)(x) > x$.
\end{lemma}
\begin{proof} Straightforward.

$\mathsf{DYNM}(x,n)(x)=\mathsf{DYNM}(2^n-1,n)(x)=I^n(x)\\
=(3(...(3(3x+1)/2)+1)/2...)+1)/2 \\
=\frac{3}{2}(\frac{3}{2}(...\frac{3}{2}(\frac{3}{2}x+\frac{1}{2})+\frac{1}{2})+...+\frac{1}{2})+\frac{1}{2}\\
=(\frac{3}{2})^nx+\frac{1}{2}((\frac{3}{2})^{n-1}+(\frac{3}{2})^{n-2}+...+1)\\
=(\frac{3}{2})^nx+\frac{1}{2}(\frac{(\frac{3}{2})^n-1}{\frac{3}{2}-1})
= (\frac{3}{2})^nx+(\frac{3}{2})^n-1 = (\frac{3}{2})^n(x+1)-1 \\
= (\frac{3}{2})^n*x+(\frac{3}{2})^n-1 >
(\frac{3}{2})^n*x > x.$ \qed
\end{proof}

\begin{lemma} \label{lemma:rdx>n}
$\forall n \in \mathbb{N}^*,$ $\exists x=2^n-1 \in \mathbb{N}^*$, if
$\mathsf{RD}[x]$ exists, then $|\mathsf{RD}[x]|>n$.
\end{lemma}
\begin{proof} Straightforward.
$\mathsf{DYNM}(x,n)(x)=\mathsf{DYNM}(2^n-1,n)(x)=I^n(x)>x.$ If
$\mathsf{RD}[x]$ exists, then $|\mathsf{RD}[x]| \geq n+1 >n$. \qed
\end{proof}

By above lemmas, we can draw a conclusion that maximal length of
reduced dynamics is nonexistent. Hence, it is \emph{impossible} to
enumerate all reduced dynamics. Indeed, original dynamics is also
held for above two conclusions.

\begin{theorem} \label{th:maxisnonexits} If Collatz conjecture is true, then
the maximal count of Collatz transformations for all positive
integers to return to 1 is nonexistent.
\end{theorem}
\begin{proof} %By Lemma \ref{lemma:rdx>n},
%$\forall n \in \mathbb{N}^*,$ $\exists x=2^n-1 \in \mathbb{N}^*$, if
%$\mathsf{RD}[x]$ exists, then $|\mathsf{RD}[x]|>n$. If Collatz
%conjecture is true, then $\forall x \in \mathbb{N}^*,$
%$\mathsf{RD}[x]$ exists. Thus, $\forall n \in \mathbb{N}^*,$
%$\exists x=2^n-1 \in \mathbb{N}^*$ such that $\mathsf{RD}[x]$ exists
%and $|\mathsf{RD}[x]|>n.$ Hence, the maximal number for all natural
%numbers cannot exists.
%
%Alternatively,

Collatz conjecture is true, hence $\forall x \in \mathbb{N}^*,$
there exists original dynamics of $x$. Let $S=\{n | \forall x \in
\mathbb{N}^*, \mathsf{DYNM}(x,n)(x)=1 \}.$ That is, $S$ is a set of
all $x \in \mathbb{N}^*$ that can return to 1 after finite times of
Collatz transformations (i.e., $n$, in terms of $(3*x+1)/2$ denoted
as $I$ and $x/2$ denoted as $O$). In other words, $\forall x \in S$,
$\exists n \in \mathbb{N}^*$ such that $\mathsf{DYNM}(x,n)(x)=1$.

Next, we will prove that $\max(S)$ cannot exist.

Suppose $\max(S)$ exists. We next construct a contradiction as
follows: Let $N_{max}=\max(S)$. We can create $y=2^{N_{max}}-1 \in
\mathbb{N}^*$. As Collatz conjecture is true, $y$ can return to 1
after finite times of Collatz transformations. Thus, $y \in S$.
However, $y$ needs at least $N_{max}+1$ times of transformations to
return to 1 due to $\mathsf{DYNM}(y,N_{max})(y)=I^{N_{max}}(y)>y$
(i.e., $\mathsf{DYNM}(y,N_{max})=I^{N_{max}}$), which contradicts
with $N_{max}=\max(S)$. \qed
\end{proof}

\begin{corollary} If Collatz conjecture is true, the set of all $x \in \mathbb{N}^*$ that
can return to 1 after finite times of Collatz transformations cannot
be enumerated.
\end{corollary}

\begin{proof}
Straightforward due to Theorem \ref{th:maxisnonexits}.

Suppose $S$ is a set of $x \in \mathbb{N}^*$ that can return to 1
after finite times of Collatz transformations in terms of
$(3*x+1)/2$ (denoted as $I$) and $x/2$ (denoted as $O$). Initially,
$S$ is empty. If $\forall x \in S$, $\exists n \in \mathbb{N}^*$
such that $\mathsf{DYNM}(x,n)(x)=1,$ then let $x$ be included into
$S$.

Let $N_{max}=\max(\{n | \forall x \in S,
\mathsf{DYNM}(x,n)(x)=1\}).$ We can create $y=2^{N_{max}}-1$. As
Collatz conjecture is true, $y$ can return to 1 after finite times
of Collatz transformations. $y \not \in S$, because $y$ needs at
least $N_{max}+1$ times of transformations to return to 1 due to
$\mathsf{DYNM}(y,N_{max})(y)=I^{N_{max}}(y)>y$ (i.e.,
$\mathsf{DYNM}(y,N_{max})=I^{N_{max}}$). That is, the set $S$ cannot
be enumerated in finite times, there always exists a positive
integer that is not in current $S$. \qed
\end{proof}

Next, we will provide an approach to prove Reduced Collatz
Conjecture.

\begin{definition} Function $CntO(\cdot)$. $CntO: s \rightarrow n$
takes as input $s \in \{I,O\}^{\geq1}$, and outputs $n\in
\mathbb{N}^*$ that is the count of ``$O$'' in $s$.
\end{definition}

\begin{example} $CntO(IIOO)=2.$ Obviously, the function name
comes from ``counting the number of $O$''.
\end{example}

We can create an one-to-one mapping between residue class and
reduced dynamics. (Indeed, in our another paper \cite{wei09ratio},
the formula to compute residue class for a given reduced dynamics is
derived. Sufficient and necessary condition for reduced dynamics is
proved. The period due to residue class is extensively explored and
proved in our another paper %\cite{wei09period}.
Those major
properties for reduced dynamics such as partition, period, and ratio
can be showed visually in our proposed graph \cite{weiispa}, which
is provided in Appendix).

A general algorithm can be proposed by Theorem \ref{th:partition},
which outputs residue class of $x$, and takes as input first $t$
transformations in (partial) reduced dynamics of $x$ (i.e.,
$s=\mathsf{DYNM}(x,t)$, $|\mathsf{RD}[x]|\geq t$).

Given specific $s \in \{I,O\}^t,$ residue class $[i]_{2^t}$ can be
determined by Algorithm \ref{alg:ct^t2i} (called D2R Algorithm) such
that $\mathsf{DYNM}(x\in [i]_{2^t},t)=s$ as follows:

\begin{algorithm}[!h] \label{alg:ct^t2i}
\caption{D2R Algorithm. Input $s \in \{I,O\}^t, t \in \mathbb{N}^*,$
$t\geq 2$, $CntO(GetS(s,1,j))<\lceil \log_21.5*CntI(GetS(s,1,j))
\rceil$ where $j=1,2,...,t-1,$ $CntO(s) \leq \lceil
\log_21.5*CntI(s) \rceil$. Output $[i]_{2^t}$, $1 \leq i \leq
2^t-1,$ $i \in [1]_2$, such that $\mathsf{DYNM}(x \in
[i]_{2^t},t)=s.$} \KwData{$s$} \KwResult{$i$}

    $i \Leftarrow 1$\;

    %\If{$t=1 \wedge f=O$} {$i\Leftarrow 0$\;}

    \For{$j=1; j \leq t-1; j++$} {

        \If{$IsMatched(GetS(s,1,j)(i), GetS(s,j+1,1))==False$} {
            $i \Leftarrow 2^j+i$\;
        }
    }
\Return $[i]_{2^t}$\;
\end{algorithm}

\begin{remark}

\item (1) If $IsMatched(GetS(s,1,j)(i), GetS(s,j+1,1))=False$,
then $i \Leftarrow 2^j+i$. It is due to Theorem \ref{th:partition}.

\item (2) Indeed, $i$ in outputs $[i]_{2^t}$ is the major computation result, as $t$ is
available in the input $s$ due to $t=|s|$.

\item (3) Obviously, above algorithm can be terminated and the
time cost is $O(t).$

\item (4) Besides, $t=1$ is trivial and can be easily included in the
algorithm as follows: If input is $I,$ then let the output of D2R
algorithm be $[1]_{2^1}$. If input is $O$, then let the output be
$[0]_{2^1}$. Consequently, Corollary \ref{cor:n1=n2} can be extended
to include $t=1$ easily.

\end{remark}

\begin{corollary} \label{cor:n1=n2}

(1) Given $s \in \{I,O\}^t, t \in \mathbb{N}^*, t \geq 2$,\\
$CntO(GetS(s,1,j))<\lceil \log_21.5*CntI(GetS(s,1,j)) \rceil$, $j=1,2,...,t-1$,\\
$CntO(s) \leq \lceil \log_21.5*CntI(s) \rceil,$\\
there exists one and only one $i \in [1,2^t-1], i \in [1]_2$ such
that\\
$\mathsf{DYNM}(x\in [i]_{2^t},t)=s.$

(2) Inversely, given $x \in [i]_{2^t}, t \in \mathbb{N}^*, t \geq 2,
i \in [1,2^t-1], i \in [1]_2$,\\
there exists one and only one $s \in \{I,O\}^t$ such that\\
$\mathsf{DYNM}(x,t)=s$, \\
$CntO(GetS(s,1,j))<\lceil \log_21.5*CntI(GetS(s,1,j))\rceil$,
$j=1,2,...,t-1$, \\
$CntO(s) \leq \lceil \log_21.5*CntI(s) \rceil.$

(3) If $N_1=\|\{s|s \in \{I,O\}^t,$ $t\in \mathbb{N}^*, t\geq 2,$ $i \in [1,2^t-1], i \in [1]_2,$ \\
$CntO(GetS(s,1,j))<\lceil \log_21.5*CntI(GetS(s,1,j))\rceil,$
$j=1,2,...,t-1,$ $CntO(s) \leq \lceil \log_21.5*CntI(s) \rceil,$
$\mathsf{DYNM}(x\in [i]_{2^t},t)=s \}\|$, and

$N_2=\|\{i|s \in \{I,O\}^t,$ $t\in \mathbb{N}^*, t\geq 2$, $i \in [1,2^t-1], i \in [1]_2,$ \\
$CntO(GetS(s,1,j))<\lceil \log_21.5*CntI(GetS(s,1,j))\rceil$, $j=1,2,...,t-1,$ \\
$CntO(s) \leq \lceil \log_21.5*CntI(s) \rceil,$ $\mathsf{DYNM}(x\in
[i]_{2^t},t)=s \}\|,$

then $N_1=N_2.$

\end{corollary}
\begin{proof}
(1) $[i]_{2^t} \Rightarrow s.$ It is due to the computation for the
first $t$ transformations of $i$, which is deterministic, and all $x
\in [i]_{2^t}$ have the same first $t$ transformations due to Lemma
\ref{lemma:IsOdd'}.
%Corollary \ref{cor:tdetermin}.

(2) $s \Rightarrow [i]_{2^t}.$ As all transformations of $[j]_{2^t},
j=1,...,2^t-1, j \in [1]_2$ can be enumerated, one and only one of
them equals $s$ due to Lemma \ref{lemma:IsOdd'}. That $j$ is $i$.
Alternatively, non-trivial algorithm outputting $i$ for given $s$ is
proposed in Algorithm \ref{alg:ct^t2i}.

(3) Due to (1), we have $N_1=N_2$. (Recall Corollary
\ref{cor:rstody1} and Corollary \ref{cor:rstody2}.) \qed
\end{proof}

\begin{remark}
\item (1) Simply speaking, above corollary states that $s$ and $[i]_{2^t}$
can be mutually determined such that $\mathsf{DYNM}(x \in
[i]_{2^t},t)=s,$ and thus the number of types are identical.

\item (2) For better understanding above corollary, we can explain or
observe above $N_1$ and $N_2$ in our proposed reduced dynamics graph
in Appendix. That is, $N_1$ is the count of paths consisting of
either ``$I$'' or ``$O$'' with length $t$ from starting integer;
$N_2$ is the count of distinct residue classes whose first $t$
transformations equal these paths.

\end{remark}

\begin{proposition} (An approach to prove Reduced Collatz Conjecture.)
Given $n \in \mathbb{N}^*,$ the number of integers $x$ ($x \in
\mathbb{N}^*$) such that $|\mathsf{RD}[x]| \leq n$ is infinite due
to periodical property of $x$. We denote this set as $S(x,n) =\{x(n)
| n,x \in \mathbb{N}^*, |\mathsf{RD}[x]| \leq n \}$. The ratio of
the number of this set over the number of all positive integers is
$R(x,n) = \frac{\|S(x,n)\|}{\|\mathbb{N}^*\|}$. It is finite and
obviously $R(x,n) \in (0,1]$. E.g., $S(x,1)=[0]_2;$ $S(x,2)=\{x(2)|
|\mathsf{RD}[x]| \leq 2\}=\{x(2)| |\mathsf{RD}[x]|=1,
|\mathsf{RD}[x]|=2\}=[0]_2 \cup [1]_4;$
$R(x,1)=\frac{\|[0]_2\|}{\|\mathbb{N}^*\|}=1/2;$
$R(x,2)=\frac{\|[0]_2 \cup [1]_4\|}{\|\mathbb{N}^*\|}=1/2+1/4=3/4.$

We hereby state that, if $R(x,i)$ goes to 1 with the growth of $n$
to infinite, then Reduced Collatz Conjecture is true (and thus
Collatz Conjecture is true). That is, if $$\lim_{n\rightarrow
+\infty} R(x,n)=1,$$ where $R(x,n)=
\frac{\|S(x,n)\|}{\|\mathbb{N}^*\|}$, then Reduced Collatz
Conjecture is true and thus Collatz Conjecture is true.
\end{proposition} \begin{proof} Straightforward. If
$\lim_{n\rightarrow +\infty} R(x,n)=1,$ then
$\|S(x,n)\|=\|\mathbb{N}^*\|,$ when $n\rightarrow +\infty.$ Thus,
$\forall x \in \mathbb{N}^*,$ $\exists \mathsf{RD}[x].$ Then,
Reduced Collatz Conjecture is true and thus Collatz Conjecture is
true. \qed \end{proof}

\section{Conclusion \label{sec:Conclusions}}
%We use $\mathsf{RD}[x]$ to denote the reduced dynamics of $x$, which
%is represented by a sequence of computation of $I$ or $O$, where $I$
%denotes $(3*x+1)/2$ and $O$ denotes $x/2$.

This paper discovered and proved that, all positive integers,
especially odd numbers, are partitioned regularly. If the first
different $f \in \{I,O\}$ transformation of $x_1,x_2$ occurs at the
$(t+1)$-th transformation, $t \in \mathbb{N}^*$, then $x_1 \equiv
x_2 \mod 2^t$ and $x_1 \equiv x_2 + 2^t \mod 2^{t+1}.$ The first $t$
transformations for all $x \in [x \mod 2^t]_{2^t}$ are identical.
The $(t+1)$-th transformation for $x$ and $x+2^t$ is distinct.
$\forall x\in [1]_2,$ if $s(x) \not < x,$ %$C2O(s) < \lceil \lambda*C2I(s)\rceil,
$s=\mathsf{DYNM}(x,t)$, $t\geq 2$, then $[i]_{2^{t}}, i \in [1]_2$
is partitioned into two halves and either half presents $I$ or $O$
in the $(t+1)$-th transformation. Otherwise, reduced dynamics is
obtained, i.e., $\mathsf{DYNM}(x,t)=\mathsf{RD}[x].$

%Finally, the union of residue classes who present reduced dynamics
%become larger with the growth of residue module, and equals all
%natural numbers asymptotically. $[ListI_U]_{m_U} \subset [1]_2
%\backslash \bigcup_{u=1}^{U-1} [ListI_u]_{m_u},$ and \\
%$\lim_{U \rightarrow +\infty}\bigcup_{u=0}^U[ListI_u]_{m_u} =
%\mathbb{N}^*.$ Therefore, Reduced Collatz Conjecture is TRUE, and
%equivalently Collatz Conjecture is TRUE.

\section*{Data Availability}

The data (Source Code in C, Computer Program Outputs) used to
support the findings of this study have been deposited in the [1-3]
repository. (1) Wei Ren, Exploring properties in reduced Collatz
dynamics, IEEE Dataport, 2018.  http://dx.doi.org/10.21227/ge9h-8a77
(2018). (2) Wei Ren, Verifying whether extremely large integer
guarantees Collatz conjecture (can return to 1 finally), IEEE
Dataport, http://dx.doi.org/10.21227/fs3z-vc10 (2018). (3) Wei Ren,
Exploring the ratio between the count of x/2 and the count of
(3*x+1)/2 in original dynamics for extremely large starting integers
asymptotically, IEEE Dataport. http://dx.doi.org/10.21227/rxx6-8322
(2018). (4) Wei Ren, Exploring the inverse mapping from a dynamics
to a residue class - inputting a reduced dynamics or partial
dynamics and outputting a residue class , IEEE Dataport.
http://dx.doi.org/10.21227/qmzw-gn71 (2018). (5) Wei Ren, Reduced
Collatz Dynamics for Integers from 3 to 999999, IEEE Dataport.
http://dx.doi.org/10.21227/hq8c-x340 (2018). (6) Wei Ren, Collatz
Automata and Compute Residue Class from Reduced Dynamics by Formula,
IEEE Dataport. http://dx.doi.org/10.21227/7z84-ms87 (2018).

\section*{Conflicts of Interest}

The author has no conflicts of interest to declare.

\section*{Acknowledgement} The research was financially supported by
the Provincial Key Research and Development Program of Hubei (No.
2020BAB105), Knowledge Innovation Program of Wuhan-Basic Research
(No. 2022010801010197), the Foundation of National Natural Science
Foundation of China (No. 61972366), and the Opening Project of
Nanchang Innovation Institute, Peking University (No. NCII2022A02).
%A preprint has %previously been published \cite{weipreprint}.

%\section*{Conflicts of Interest}
%
%The author has no conflicts of interest to declare.
%
\clearpage
\section*{Appendix: Dynamics Graph with Partition Labels}

\begin{figure}[!ht]
\begin{centering}
\includegraphics[width=14cm]{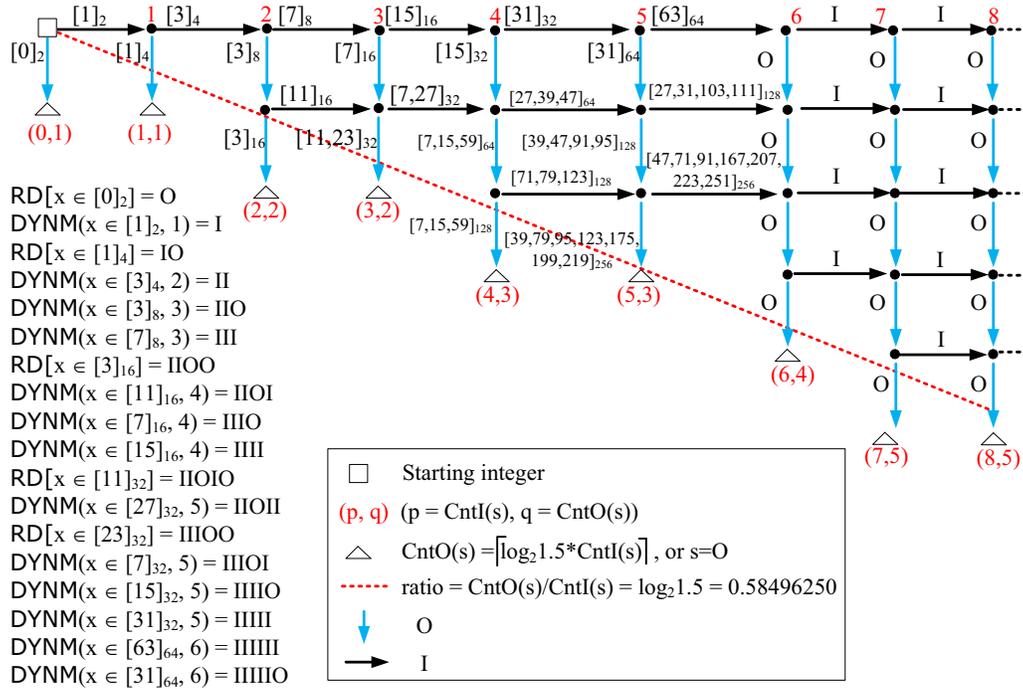}
\par\end{centering}
\caption{Reduced dynamics graph with partition labels. ``$\square$''
represents starting integer. ``$\triangle$'' represents transformed
integer that is less than starting integer. Any ``$\triangle$'' is
below ratio line, whose slope is $\lambda=\log_21.5$. A reduced
dynamics is a path consisting of edges in terms of ``$I$'' or
``$O$'' from ``$\square$'' to ``$\triangle$''. If and only if
$|s|=1, s=O$ or $CntO(s)=\lceil \lambda*CntI(s)\rceil$ and
$CntO(s')<\lceil \lambda*CntI(s') \rceil, s'=GetS(s,1,j),
j=1,2,...,|s|-1, |s| \geq 2$, then $s$ is a reduced dynamics, i.e.,
$s(x)<x,$ $GetS(s,1,j)(x) \not<x,$ $j=1,2,...,|s|-1.$ All reduced
dynamics intersects the ratio line in the last edge, which means
transformed integer is the first one that is less than corresponding
starting integer. \label{fig:allcodebranchwithx.eps}}
\end{figure}

Reduced dynamics graph \cite{weiispa} with partition labels can be
obtained by adding each branch of paths a partition of $x \in
\mathbb{N}^*$ (see Fig.\ref{fig:allcodebranchwithx.eps}), which
provides a visualization and smooth understanding for our proofs of Partition Theorem. %(Section
%\ref{sec:partition}.

%The graph presents information as follows:
\begin{remark}
\item (1) Either partition of current $x$ presents the same previous
dynamics, but will present either edge ($I$ or $O$) in the next
transformation, iteratively.

More specifically, either half of $x \in [i]_{2^t}$ (namely, $x \in
[i]_{2^{t+1}}$ or $x \in [i+2^t]_{2^{t+1}}$) has either
$\mathsf{DYNM}(x,t+1)=\textsf{DYNM}(x,t)\|I$ or
$\mathsf{DYNM}(x,t+1)=\mathsf{DYNM}(x,t)\|O$. Besides, $\forall x
\in [i]_{2^t},\mathsf{DYNM}(x,t)$ are identical.

For example, $\mathsf{DYNM}(x \in [3]_8,3)=IIO$. Whether the next
transformation is $I$ or $O$ depends on either half partition of $x
\in [3]_8$. That is, if $x \in [3]_{16}$, then $IIO(x) \in [0]_2$
and the next one is $O$; If $x \in [11]_{16},$ then $IIO(x) \in
[1]_2$ and the next one is $I$. Obviously, $[3]_8 =
[3]_{16}\cup[11]_{16}$ and $[3]_{16}\cap [11]_{16}=\emptyset$.

\item (2) Let $\mathsf{DYNM}(x,L)=s \in \{I,O\}^L, L \in \mathbb{N}^*.$ If
$CntO(s)=\lceil \lambda*CntI(s)\rceil$ (or $s=O$), then final
transformed integers that are denoted as ``$\triangle$'' is avaiable
(i.e., $\mathsf{DYNM}(x,L) = \mathsf{RD}[x]$. In other words,
$s(x)<x,$ $GetS(s,1,j)(x) \not<x,j=1,2,...,L-1$).

\item (3) Each edge (namely, $I$ or $O$) is labeled with one residue class
or a union of multiple residue classes, whose last transformation is
this edge ($I$ or $O$), and $x$ in this residue class presents the
dynamics represented by the path from ``$\square$'' (representing
starting integer) to this edge.

For example, $[11,23]_{32}$ is labeled due to $\mathsf{DYNM}(x \in
[11,23]_{32},5)=s_{11},s_{23} \in \{I,O\}^5$ where
$GetS(s_{11},5,1)=O, GetS(s_{23},5,1)=O.$ More specifically,
$\mathsf{DYNM}(x \in [11]_{32},5)=s_{11}=IIOIO$, $\mathsf{DYNM}(x
\in [23]_{32},5)=s_{23}=IIIOO$.

More specifically, if $[i]_{2^t} (or [i_1,i_2,...,i_p]_{2^t}, p \in
\mathbb{N}^*$) is (are) residue class (classes) as a label located
at the last edge of path $s_{i} \in \{I,O\}^t$ (or paths
$s_{i_1},s_{i_2},...,s_{i_p}$), then $\mathsf{DYNM}(x \in
[i]_{2^t},t)=s_i$ (or $\mathsf{DYNM}(x \in
[i_j]_{2^t},t)=s_{i_j},j=1,2,...,p$).

\end{remark}

For reduced dynamics, the ratio - the count of $x/2$ (or ``$O$'')
over the count of $(3*x+1)/2$ (or ``$I$'') - is larger than
$\log_21.5$ (denoted as $\lambda$). That is, when and only when the
count of ``$O$'' is larger than the count of ``$I$'' times
$\lambda$, current transformed integer will be less than the
starting integer (i.e., the reduced dynamics will be available).
Indeed, we prove following sufficient and necessary condition in
which any $s \in \{I,O\}^{\geq1}$ is a reduced dynamics formally in
another paper \cite{wei09ratio}.

\begin{corollary} \label{th:codetheoremsufficiency} (\textbf{Form
Corollary.}) $s \in \{I,O\}^{\geq1}$ is a reduced dynamics, if and
only if

(1) $|s|=1, s=O;$ Or,

(2) $|s|\geq 2,$
\begin{equation} \left\{ \begin{aligned}
CntO(s) = \lceil \lambda*CntI(s) \rceil, && \\
CntO(s')< \lceil \lambda*CntI(s') \rceil, && s'=GetS(s,1,i),
i=1,2,...,|s|-1.
\end{aligned} \right.
\end{equation}
\end{corollary}


\begin{thebibliography}{00}

%% \bibitem{label}
%% Text of bibliographic item

%\bibitem{Livio11}
%L. Colussi. \emph{The convergence classes of Collatz function}
%\newblock{Theoretical Computer Science}, vol. 412, no. 39, pp.
%5409-5419, 2011
%
%\bibitem{Hew16}
%P.C. Hew. \emph{Working in binary protects the repetends of 1/3h:
%Comment on Colussi's 'The convergence classes of Collatz function'}.
%\newblock{Theoretical Computer Science}, vol. 618, pp. 135-141, 2016

%\bibitem{Guy1983Don}
%R.K. Guy, \emph{Don't try to solve these problems},
%\newblock{Computers and Mathematics with Applications}, vol. 90, no.
%1, pp. 35-41, 1983.



\bibitem{marc06} Marc Chamberland, \emph{An Update on the 3x+1
Problem} \newblock{Butlleti de la Societat Catalana de
Matematiques}, 18, pp.19-45,
https://chamberland.math.grinnell.edu/papers/3x\_survey\_eng.pdf

\bibitem{lagarisa10} Jeffrey C. Lagarias, \emph{The 3x+1 Problem: an
Overview}, \newblock{The Ultimate Challenge: The 3x + 1 Problem},
Edited by Jeffrey C. Lagarias. American Mathematical Society,
Providence, RI, 2010, pp. 3-29,
https://doi.org/10.48550/arXiv.2111.02635

\bibitem{lagarisa11} Jeffrey C. Lagarias, \emph{The 3x + 1 Problem:
An Annotated Bibliography (1963šC1999)}, \newblock{January 1, 2011
version)} https://doi.org/10.48550/arXiv.math/0309224

\bibitem{lagarisa12} Jeffrey C. Lagarias, \emph{The 3x+1 Problem: An
Annotated Bibliography, II (2000-2009)}, \newblock{January 10, 2012
version)} https://doi.org/10.48550/arXiv.math/0608208

\bibitem{UpperboundRecord1} Tomas Oliveira e Silva, \emph{Maximum
excursion and stopping time record-holders for the 3x+1 problem:
computational results}, \newblock{\em Mathematics of Computation},
vol. 68, no. 225, pp. 371-384, 1999.

\bibitem{UpperboundRecord2} Tomas Oliveira e Silva, \emph{Empirical
Verification of the 3x+1 and Related Conjectures}. In \newblock{\em
The Ultimate Challenge: The 3x+1 Problem}, (book edited by Jeffrey
C. Lagarias), pp. 189-207, AMS, 2010. American Mathematical Society,
2010.

\bibitem{barina21} David Barina, \emph{Convergence verification of
the Collatz problem}, \newblock{The Journal of Supercomputing}, 77,
2681-2688, 2021. https://doi.org/10.1007/s11227-020-03368-x

\bibitem{eric} Eric Roosendaal, \emph{On the 3x + 1 problem},
\newblock{http://www.ericr.nl/wondrous/index.html}

\bibitem{weiuic} Wei Ren, Simin Li, Ruiyang Xiao and Wei Bi,
\emph{Collatz Conjecture for $2^{100000}-1$ is True - Algorithms for
Verifying Extremely Large Numbers}, \newblock{Proc. of IEEE UIC},
Oct. 2018, Guangzhou, China, 411-416, 2018


\bibitem{Livio11}
Livio Colussi. \emph{The convergence classes of Collatz function}
\newblock{Theoretical Computer Science}, vol. 412, no. 39, pp.
5409-5419, 2011

\bibitem{Hew16}
Patrick Chisan Hew. \emph{Working in binary protects the repetends
of 1/3h: Comment on Colussi's 'The convergence classes of Collatz
function'}. \newblock{Theoretical Computer Science}, vol. 618, pp.
135-141, 2016

\bibitem{Leavens} G.T. Leavens and M. Vermeulen, \emph{3x+1 Search
Programs}, \newblock{Computers and Mathematics with Applications},
vol. 24, no. 11, pp. 79-99, 1992.

\bibitem{Crandall} R.E. Crandall, \emph{On the $``3x+1''$ problem},
\newblock{Mathematics of Computation}, vol. 32, no. 144, pp.
1281-1292, 1978.

\bibitem{Krasikov} I. Krasikov and J.C. Lagarias, \emph{Bounds for
the 3x + 1 problem using difference inequalities}, \newblock{Acta
Arithmetica}, vol. 109, no. 3, pp. 237-258, 2003.

\bibitem{weijm} Wei Ren, \emph{A New Approach on Proving Collatz
Conjecture}, \newblock{Journal of Mathsmatics}, Hindawi, April 2019,
ID 6129836, https://doi.org/10.1155/2019/6129836.

\bibitem{weiispa} Wei Ren, \emph{Ratio and Partition are Revealed in
Proposed Graph on Reduced Collatz Dynamics}, \newblock{Proc. of 2019
IEEE Intl Conf on Parallel \& Distributed Processing with
Applications (ISPA)}, pp. 474-483, 16-28 Dec. 2019, Ximen, China

%\bibitem{weihpcc}
%Wei Ren, Ruiyang Xiao, \emph{How to Fast Verify Collatz Conjecture
%by Automata}, \newblock{Proc. of %IEEE 21st International Conference
%%on High Performance Computing and Communications; IEEE 17th
%%International Conference on Smart City; IEEE 5th International
%%Conference on Data Science and Systems (HPCC/SmartCity/DSS)
%IEEE HPCC}, pp. 2720-2729, 10-12 Aug. 2019, Zhangjiajie, China

%\bibitem{weidata}
%Wei Ren, \emph{Reduced Collatz Dynamics Data Reveals Properties for
%the Future Proof of Collatz Conjecture},
%\newblock{Data}, MDPI, 2019, 4, 89. %doi:10.3390/data4020089,
%%https://www.mdpi.com/2306-5729/4/2/89/pdf.

%\bibitem{wei09period}
%Wei Ren, \emph{Reduced Collatz Dynamics is Periodical and the Period
%Equals 2 to the Power of the Count of x/2},
%\newblock{Computational and Mathematical Methods}, Wiley and Hindawi, Submitted, 2022.

\bibitem{wei09ratio} Wei Ren, \emph{A Reduced Collatz Dynamics Maps
to a Residue Class, and its Count of x/2 over Count of 3*x+1 is
larger than ln3/ln2}, \newblock{International Journal of Mathematics
and Mathematical Sciences}, Hindawi, Volume 2020, Article ID
5946759, https://doi.org/10.1155/2020/5946759, 2020.

%\bibitem{weipreprint}
%Wei Ren, \emph{Collatz Dynamics is Partitioned by Residue Class
%Regularly}, \newblock{Preprint},
%https://www.preprints.org/manuscript/202209.0155/v1, 2022


%\bibitem{wei09partition}
%Wei Ren, \emph{Collatz Dynamics is Partitioned by Residue Class
%Regularly},
%\newblock{Information and Computation}, Elsevier, Submitted, 2019.


%\bibitem{Livio11}
%Livio Colussi. \emph{The convergence classes of Collatz function}
%\newblock{Theoretical Computer Science}, vol. 412, no. 39, pp.
%5409-5419, 2011
%
%\bibitem{Hew16}
%Patrick Chisan Hew. \emph{Working in binary protects the repetends
%of 1/3h: Comment on Colussi's 'The convergence classes of Collatz
%function'}. \newblock{Theoretical Computer Science}, vol. 618, pp.
%135-141, 2016

%\bibitem{chap1-5} Wei Ren.
%\newblock Induction and Code for Collatz Conjecture or 3x+1 Problem.
%\newblock{\em Submitted to Journal of Algebra, Elsevier. ArXiv.org},
%
%\bibitem{chap6} Wei Ren.
%\newblock The Quantitative Law of Dynamics in 3x+1 Conjecture - $D \geq \lceil
%\log_23*U \rceil$.
%\newblock{\em Submitted to Results in Mathematics, Springer. ArXiv.org},
%
%\bibitem{chap7} Wei Ren.
%\newblock Dynamics is Periodical and Period is $2^D$ in 3x+1 Conjecture - $CODE(x+2^D)=CODE(x)$.
%\newblock{\em Submitted to Discrete Mathematics, Elsevier. ArXiv.org},

%\bibitem{chap9} Wei Ren.
%\newblock Reversibility of $\mathsf{RD}[x]$ - Mapping Dynamics to Residue Classes in 3x+1 Conjecture
%\newblock{\em Submitted to Discrete Applied Mathesmatics, Elsevier. ArXiv.org},

%\bibitem{chap11} Wei Ren.
%\newblock Collatz Conjecture for $2^{100000}-1$ is True - Algorithms for Verifying Extremely Large Numbers.
%\newblock{\em Submitted to Computing, Springer. ArXiv.org},

%\bibitem{chap12} Wei Ren.
%\newblock Output Dynamics of Collatz Conjecture by Automata.
%\newblock{\em Submitted to Theoretical Computer Science, Elsevier. ArXiv.org},

\end{thebibliography}
\end{document}